\documentclass[a4paper,12pt]{amsart}
\usepackage{latexsym,amsmath,amssymb,amsthm}
\newtheorem{thm}{Theorem}[section]
\newtheorem{lemma}[thm]{Lemma}

\newtheorem{rem}[thm]{Remark} 

\newtheorem{define}[thm]{Definition}
\newtheorem{prop}[thm]{Proposition}

\numberwithin{equation}{section}

\def\a{\alpha}
\def\b{\beta}
\def\e{\varepsilon}

\def\H{\mathcal{H}}

\def\F{\mathcal{F}}
\def\o{\omega}

\def\p{\partial}

\def\R{\mathbb R}
\def\T{\mathbb T}

\def\N{\mathbb N}

\input{amssym.def}

\title[A second derivative H\"{o}lder estimate for weak MCF]{A second  derivative H\"{o}lder estimate \\ for weak mean curvature flow}
\author[Y. Tonegawa]{Yoshihiro Tonegawa}

\address{Department of Mathematics, Hokkaido University, Sapporo 060-0810 Japan.}
\email{tonegawa@math.sci.hokudai.ac.jp}

\date{}


\thanks{Partially supported by JSPS Grant-in-aid for scientific research (B) $\#$21340033, (S) $\#$21224001
and challenging exploratory research $\#$23654057.}

\begin{document}
\setlength\parskip{5pt}
\begin{abstract}
We give a proof that Brakke's mean curvature flow under the unit density assumption
is smooth almost everywhere 
in space-time. More generally, if the velocity is equal in a weak sense to its 
mean curvature plus some given $\a$-H\"{o}lder continuous vector field, then we show $C^{2,\a}$
regularity almost everywhere. 
\end{abstract}

\maketitle
\makeatletter
\@addtoreset{equation}{section}                  
\renewcommand{\theequation}{\thesection.\@arabic\c@equation}
\makeatother

\section{Introduction}
A family $\{M_t\}_{t\geq 0}$ of $k$-dimensional surfaces in $\R^n$ is called the mean curvature
flow (hereafter abbreviated MCF) if the velocity of $M_t$ is equal to its mean curvature at each point and
time. The MCF has been the subject of intensive research since 1980's due to its 
importance in the analytic and geometric context as well as for various applications to physical and information
sciences such as image
processing and metallurgy. The most pertinent aspect of MCF to the present
paper is the fact that the MCF is the natural gradient flow of the $k$-dimensional surface area and hence is 
equipped with uniquely rich variational structures.  In his seminal work \cite{Brakke}, Brakke took the 
advantage to define and study his version of MCF, so called Brakke's MCF 
(or we may call `weak MCF' to include more general flows), using the notion
of varifold \cite{Allard} in geometric measure theory. 
More precisely, given any $k$-dimensional integral varifold $V_0$, 
which may be considered as a generalized $k$-dimensional surface with possible singularities,  
Brakke proved the existence of a family of varifolds $\{V_t\}_{t\geq 0}$ each of which satisfies the MCF
equation taking the advantage of its variational characterization. Under the further assumption that the
density function is 1 almost everywhere in time and space, 
Brakke also claimed that the MCF is smooth almost everywhere
and that it satisfies the MCF equation in the classical sense. The
proof of regularity theorem contains remarkable new insights such as `clearing-out',
`popping soap film' and `cylindrical growth rates', to name a few.  On the other hand it is 
technically involved and  some part, in particular the graphical approximations
of the support of moving varifolds \cite[Sec. 6.9, `Flattening out']{Brakke},  is particularly difficult to follow.
Later a local regularity theorem for special but very useful case was obtained
by White \cite{White1} which is sufficient for many applications of interest while
it does not replace Brakke's claims in full. Recently Kasai and the author \cite{Kasai-Tonegawa}
gave a new proof for Brakke's regularity theorem up to $C^{1,\varsigma}$ for general weak MCF where
the velocity can be equal to the mean curvature plus any given ambient vector field in a suitable integrability class.
Note that $C^{1,\varsigma}$ here means $C^{1,\varsigma}$ in the space variables and $C^{\frac{1+\varsigma}{2}}$
in the time variable, which are the usual regularity features of parabolic problems (in the following $C^{2,\alpha}$
should be understood in the similar manner). 
The additional different aspect of \cite{Kasai-Tonegawa} from Brakke's result is that it
is a natural parabolic generalization of Allard's regularity theorem for varifold \cite{Allard}
since the time-independent case of \cite{Kasai-Tonegawa} reduces essentially to Allard's theorem. 
The new decisive input to the proof of \cite{Kasai-Tonegawa} is Huisken's monotonicity formula for MCF \cite{Huisken} 
and its variants which were not known at the time that Brakke obtained his result. 

The purpose of the present paper is to extend the regularity result from $C^{1,\varsigma}$ to $C^{2,\alpha}$ for
Brakke's MCF and more generally for weak MCF with $C^{\alpha}$ transport term. In the case of Brakke's MCF,
that is, the case that the transport term is identically equal to $0$, 
$C^{2,\alpha}$ regularity implies $C^{\infty}$ almost everywhere by the standard linear parabolic regularity theory. 
This proves Brakke's original
claim of almost everywhere $C^{\infty}$ regularity for his MCF. We noted in \cite{Kasai-Tonegawa} that
there is an essential gap in \cite{Brakke} for the step of obtaining $C^2$ regularity from $C^{1,\varsigma}$
(see \cite[Sec. 10.1]{Kasai-Tonegawa}). 
The present paper thus remedies the situation and proves that Brakke's claim was correct after all. 
Just to avoid a possible 
confusion for the reader, we should point out that $C^{1,\varsigma}$ regularity of \cite{Kasai-Tonegawa} does not
imply $C^{2,\alpha}$ simply by the standard linear parabolic regularity theory. This is because 
Brakke's formulation only gives variational inequality even with $C^{1,\varsigma}$
estimates, and not equality, thus requiring further
nonlinear analysis different from simple applications of linear theory. 

We briefly describe the method of proof. We first recall the method in \cite{Kasai-Tonegawa} for the close relevance. 
For obtaining $C^{1,\varsigma}$ regularity there, 
we used the so called blow-up argument. The essence of this argument is that one measures the deviation of
moving varifolds from some graph of affine function and proves that the deviation is closely approximated
by some graph of solution for the heat equation. If this can be established, then one has a way to take a much better 
affine function approximation to the moving varifolds in a smaller region. The iteration procedure then gives
$C^{1,\varsigma}$ estimate of the graph representing the support of moving varifolds. The strategy of the 
present paper is to measure the deviation of moving varifolds
from some graph of polynomial function which is quadratic (respectively, linear) in the space (respectively, time) variables and which
satisfies the heat equation, and to prove that the small deviation is closely approximated by some graph of solution for
the heat equation. Then one can find a much better approximation by a similar polynomial function in a smaller
region, and the iteration argument gives $C^{2,\alpha}$ estimates. 
The procedure takes advantage of $C^{1,\varsigma}$ estimate of \cite{Kasai-Tonegawa}, another version of 
$L^2$-$L^{\infty}$ type estimate different from \cite[Sec. 6.2]{Kasai-Tonegawa}, blow-up argument and 
it is similar to $C^{1,\varsigma}$ estimate in spirit. Since we already know that the support of moving 
varifolds is a $C^{1,\varsigma}$ graph, we need no Lipschitz graph approximation as was done in 
\cite{Kasai-Tonegawa}. Thus the proof is less technical in that respect but more so due to the higher
order approximations. 

There have been numerous works \cite{AWT,CGG,ES1,ES2,LS} which show the existence of generalized MCF
past singularities and global in time, 
and we see a significant advance of understandings for the special but important subclass of
mean convex hypersurfaces \cite{Metzger,White2,White3}. Numerous works which have even more direct 
relations to Brakke's MCF are singular perturbation limit problems such as the Allen-Cahn equation \cite{Ilmanen,
LST}
and the parabolic Ginzburg-Landau equation \cite{Ambrosio-Soner,BOS,Jerrard-Soner,Lin}. See \cite{Kasai-Tonegawa}
for further discussion. 
We cite \cite{Ecker} as one of the best references for Brakke's MCF. 

The organization of the paper is as follows. Section 2 contains basic definitions and notations. Section 3 
describes the assumptions and main results of the paper. Section 4 gives the supremum 
and Dirichlet energy estimate for the difference of heights between MCF graph and a 
certain quadratic function in terms of their $L^2$-norm
in a larger domain. The estimate is essentially used in the subsequent Section 5, where a blow-up
argument shows a decay estimate necessary for $C^{2,\alpha}$ estimate. Section 6 concludes the
proof of $C^{2,\a}$ estimate and Section 7 describes the application to MCF in submanifold. The last
Section 8 contains some technical estimates concerning the change of second derivatives under 
orthogonal rotations. 
 
\section{Preliminaries}
Even though the content of this section is more or less identical 
to \cite[Sec. 2]{Kasai-Tonegawa}, 
we include this section with a few changes for the reader's convenience. 
\subsection{Basic notations}
Throughout this paper, $k$ and $n$ will be positive integers with $0<k<n$. We often identify 
$\R^k$ with $\R^k\times\{0\}\subset \R^n$. Let $\N$ be the natural number
and $\R^+:=\{x\geq 0\}$. For $0<r<\infty$ and $a\in \R^n$ (or $\R^k$) let
\begin{equation*}
B_r(a):=\{x\in \R^n\,\,:\,\,|x-a|<r\},\,\,\,B_r^k(a):=\{x\in \R^k\,\,:\,\,|x-a|<r\}
\end{equation*}
and when $a=0$ let $B_r:=B_r(0)$ and $B_r^k:=B_r^k(0)$. We denote by $\H^k$ the $k$-dimensional
Hausdorff measure on $\R^n$. The restriction of $\H^k$ to a set $A$ is denoted by $\H^k\lfloor_A$.
Set $\omega_k:=\H^k(B_1^k)$. 
For an open subset $U\subset {\mathbb R}^n$ let $C_c(U)$ be the set of all compactly 
supported continuous functions on $U$ and let $C_c(U;{\mathbb R}^n)$ be 
the set of all compactly supported, continuous vector fields. The upper subscript
of $C_c^l(U)$ and $C_c^l(U;{\mathbb R}^n)$ indicates 
continuous $l$-th order differentiability. For $g\in C^1(U;{\mathbb R}^n)$,
we regard $\nabla g(x)$ as an element of ${\rm Hom}({\mathbb R}^n,{\mathbb R}^n)$.
Similarly for $g\in C^1(U)$, we regard the Hessian matrix $\nabla^2 g(x)$ as an 
element of ${\rm Hom}({\mathbb R}^n,{\mathbb R}^n)$. $\nabla$ always indicates
differentiation with respect to the space variables $x$, and not with respect to the time variable $t$. 

For any Radon measure $\mu$ on 
${\mathbb R}^n$ and $\phi\in C_c({\mathbb R}^n)$ we often write 
$\mu(\phi)$ for $\int_{{\mathbb R}^n}\phi\, d\mu$. Let ${\rm spt}\, \mu$ be
the support of $\mu$, i.e., $x\in {\rm spt}\, \mu$ if
$\mu(B_r(x))>0$ for all $r>0$. 
Let $\Theta^k(\mu,x)$ be the $k$-dimensional density of $\mu$ at $x$, i.e., 
$\lim_{r\rightarrow 0}\mu(B_r(x))/(\o_k r^k)$, when the limit exists. 
For $\mu$ a.e$.$ defined 
function $u$, and $1\leq p\leq \infty$, $u\in L^p(\mu)$ means $\left(\int |u|^p\, d\mu\right)^{1/p}<\infty$. 

For $-\infty<t<s<\infty$ and $x,\, y\in {\mathbb R}^n$, define
\begin{equation}
\rho_{(y,s)}(x,t):=\frac{1}{(4\pi(s-t))^{k/2}}\exp\left(-
\frac{|x-y|^2}{4(s-t)}\right).
\label{hkneldef}
\end{equation}
$\rho_{(y,s)}$ is the $k$-dimensional backward heat kernel.

\subsection{The Grassmann manifold and varifolds}
Let ${\bf G}(n,k)$ be the space of $k$-dimensional subspaces of ${\mathbb R}^n$
and let ${\bf A}(n,k)$ be the space of $k$-dimensional affine planes of ${\mathbb R}^n$.
For $S\in {\bf G}(n,k)$, we identify $S$ with the corresponding orthogonal 
projection of ${\mathbb R}^n$ onto $S$. Let $S^{\perp}\in {\bf G}(n,n-k)$ be the orthogonal
complement of $S$. For two  elements $A$
and $B$ of ${\rm Hom}\, ({\mathbb R}^n,{\mathbb R}^n)$, define a scalar product $A\cdot B:={\rm trace}\, (A^*\circ B)$ where $A^*$ is the transpose of $A$ and $\circ$ indicates the usual composition. The identity of ${\rm Hom}\, ({\mathbb R}^n,{\mathbb R}^n)$ is denoted by $I$. Let $a\otimes
b\in {\rm Hom}\, ({\mathbb R}^n,{\mathbb R}^n)$ be the tensor product of $a,\, b\in {\mathbb R}^n$. For $A\in {\rm Hom}\, ({\mathbb R}^n,{\mathbb R}^n)$ define
\begin{equation*}
|A|:=\sqrt{A\cdot A},\hspace{1cm}
\|A\|:=\sup\{|A(x)|\, :\, x\in {\mathbb R}^n,\, |x|=1\}.
\end{equation*}
For $T\in {\bf G}(n,k)$, $a\in {\mathbb R}^n$ and $0<r<\infty$ we define the cylinder
\begin{equation*}
C(T,a,r):=\{x\in {\mathbb R}^n\, :\, |T(x-a)|<r\},\,\, C(T,r):=C(T,0,r).
\end{equation*}
We recall some notions related to varifold and refer to \cite{Allard,Simon} for more
details. For any open set $U\subset {\mathbb R}^n$, define $G_k(U):=U\times {\bf G}(n,k)$. 
A general $k$-varifold in $U$ is a Radon measure on 
$G_k(U)$. Set of all general $k$-varifolds in $U$ is denoted by ${\bf V}_k(U)$. 
For $V\in {\bf V}_k(U)$, let $\|V\|$ be the mass measure of $V$, 
namely, 
\begin{equation*}
\|V\|(\phi):=\int_{G_k(U)}\phi(x)\, dV(x,S),\,\,\,\forall \phi\in C_c(U).
\end{equation*}
Given any ${\mathcal H}^k$ measurable 
countably $k$-rectifiable set $M\subset U$ with
locally finite ${\mathcal H}^k$ measure, there
is a natural $k$-varifold $|M|\in {\bf V}_k(U)$ defined by
\begin{equation*}
|M|(\phi):=\int_{M}\phi(x,{\rm Tan}_x M)\, d{\mathcal H}^k(x),\,\,\,\forall \phi\in C_c(G_k(U)),
\end{equation*}
where ${\rm Tan}_x M\in {\bf G}(n,k)$ is the approximate tangent space which exists
${\mathcal H}^k$ a.e$.$ on $M$. In this case, $\||M|\|={\mathcal H}^k\lfloor_M$.
We say $V\in {\bf V}_k(U)$ is integral if 
\begin{equation*}
V(\phi)=\int_{M}\phi(x,{\rm Tan}_x M)\theta(x)\, d{\mathcal H}^k(x),\,\,\,\forall \phi\in C_c(G_k(U)),
\end{equation*}
with some ${\mathcal H}^k$ measurable 
countably $k$-rectifiable set $M\subset U$ and ${\mathcal H}^k$
a.e$.$ integer-valued integrable 
function $\theta$ defined on $M$. Note that for such varifold, $\Theta^k(\|V\|,x)
=\theta(x)\in {\mathbb N}$, ${\mathcal H}^k$ a.e$.$ on $M$. Set of all 
integral $k$-varifolds in $U$ is denoted by ${\bf IV}_k(U)$.
We say $V$ is a unit density 
$k$-varifold if $V$ is integral and $\theta=1$ a.e$.$ on $M$, that is, $V=|M|$. 
When $V$ is integral, we often
write
$\int_{U}(g(x))^{\perp}\, d\|V\|(x)$ for $\int_{G_k(U)}
S^{\perp}(g(x))\, dV(x,S)$, for example, since there should be no ambiguity. 

\subsection{First variation and generalized mean curvature}
For $V\in {\bf V}_k(U)$ let $\delta V$ be the first variation of $V$, namely,
\begin{equation*}
\delta V(g):=\int_{G_k(U)}\nabla g(x)\cdot S\, dV(x,S)
\end{equation*}
for $g\in C_c^1(U;{\mathbb R}^n)$. 
Let $\|\delta V\|$ be the total variation when it exists, and if 
$\|\delta V\|$ is absolutely continuous with respect to $\|V\|$, we have
for some $\|V\|$ measurable vector field $h(V,\cdot)$
\begin{equation}
\delta V(g)=-\int_{U}g(x)\cdot h(V,x)\, d\|V\|(x).
\label{first}
\end{equation}
The vector field $h(V,\cdot)$ is called the generalized mean curvature of $V$. 
We say $V$ is stationary if $h(V,\cdot)=0$, $\|V\|$ a.e$.$ in $U$, or equivalently, 
$\delta V(g)=0$ for all $g\in C^1_c(U;{\mathbb R}^n)$. 
For any $V\in {\bf IV}_k(U)$ with integrable $h(V,\cdot)$, Brakke's 
perpendicularity theorem of generalized mean curvature
\cite[Chapter 5]{Brakke} says that we have
\begin{equation}
\int_{U} (g(x))^{\perp}\cdot h(V,x)\, d\|V\|(x)=\int_{U} g(x)\cdot h(V,x)\, d\|V\|(x)
\label{perpH}
\end{equation}
for all $g\in C_c(U;{\mathbb R}^n)$.

\subsection{The right-hand side of MCF equation}
For any $V\in {\bf V}_k(U)$, $u\in L^2(\|V\|)$
and $\phi\in C^1_c(U;{\mathbb R}^+)$, define
\begin{equation}
{\mathcal B}(V,u,\phi):=\int_{U} (-\phi(x)h(V,x)+\nabla\phi(x))\cdot(h(V,x)+
(u(x))^{\perp})\, d\|V\|(x)
\label{Bdef}
\end{equation}
when $V\in {\bf IV}_k(U)$, $\|\delta V\|$ is locally finite and absolutely continuous
with respect to $\|V\|$, and $h(V,\cdot)\in L^2(\|V\|)$. Otherwise we define
${\mathcal B}(V,u,\phi)=-\infty$. 
Formally, if a family of smooth $k$-dimensional surfaces $\{M_t\}$ moves
by the velocity equal to the mean curvature plus smooth $u$, then, one
can check that $V_t=|M_t|$ satisfies
\begin{equation}
\frac{d}{dt}\|V_t\|(\phi)\leq {\mathcal B}(V_t,u(\cdot,t),\phi),\,\,\,\forall\phi\in C_c^1(U;{\mathbb R}^+).
\label{formal}
\end{equation}
In fact, \eqref{formal}
holds with equality. Conversely, if \eqref{formal} is satisfied, then one can 
prove that the velocity is equal to the mean curvature plus $u$. If we allow the
time-varying test function $\phi\in C^1(U\times (0,\infty);\R^+)$ with $\phi(\cdot, t)\in C^1_c(U)$, 
one can check that we also have
\begin{equation}
\frac{d}{dt}\|V_t\|(\phi(\cdot,t))\leq {\mathcal B}(V_t,u(\cdot,t),\phi(\cdot,t))+\int \frac{\p\phi}{\p t}(\cdot,t)\, d\|V_t\|.
\label{formal2}
\end{equation}
This inequality \eqref{formal2} motivates the integral formulation of the motion law \eqref{main} below. 

\subsection{Notations related to norms}
For $0<\alpha<1$, $U\subset\R^n$, $-\infty<t_1<t_2<\infty$ and for any function $f: U\times(t_1,t_2)\rightarrow\R$ 
we define  the $\alpha$-H\"{o}lder semi-norm
\begin{equation*}
[ f ]_{\alpha}:=\sup_{x,\, y\in U,\, t_1<s_1<s_2<t_2} \frac{|f(x,s_1)-f(y,s_2)|}{\max\{|x-y|^{\alpha},\,|s_2-s_1|^{\alpha/2}\} }.
\end{equation*}
Though we do not write out the domain of $f$ for the notation, we always implicitly assume that the supremum is taken over the 
domain. We similarly define $[\cdot ]_{\alpha}$ for vector-valued functions and matrix-valued functions. 
For $f:U\times(t_1,t_2)\rightarrow \R$ (or $\R^m$) we also take the liberty of denoting
\begin{equation*}
\|f\|_{0}:=\sup_{x\in U,\, t\in (t_1,t_2)}|f(x,t)|
\end{equation*}
since we use $\sup$ norm quite often. Whenever it is important for clarity to specify the domain of definition, we 
write out the information. We also define the $\alpha$-H\"{o}lder norm 
\begin{equation*}
\|f\|_{\alpha}(=\|f\|_{C^{\alpha}(U\times(t_1,t_2))}):=\|f\|_{0}+[ f ]_{\alpha}.
\end{equation*}
We note that we have some occasions to define $\|f\|_{\alpha}$ differently so that it becomes scale invariant.  
This will be specified individually.  
\label{held}

\section{Main results}
\subsection{Assumptions}
For an open set $U\subset \R^n$ and $0<\Lambda\leq \infty$ suppose that we have a family 
of $k$-varifolds $\{V_t\}_{0\leq t<\Lambda}$ and a family of $n$-vector valued functions $\{u(\cdot,t)\}_{0\leq t<\Lambda}$
both on $U$ satisfying the followings.
\newline
{\bf (B1)} For a.e$.$ $t\in [0,\Lambda)$, $V_t$ is a unit density $k$-varifold.
\newline
{\bf (B2)} For $\tilde{U}\subset\subset U$ and $(t_1,t_2)\subset \subset (0,\Lambda)$, 
\begin{equation}
\sup_{t_1\leq t\leq t_2}\|V_t\|(\tilde{U})<\infty.
\label{masft}
\end{equation}
{\bf (B3)} For $0<\alpha<1$ assume that $u$ is locally $\alpha$-H\"{o}lder continuous, namely for any
$\tilde{U}\subset\subset U$ and $(t_1,t_2)\subset\subset (0,\Lambda)$,
\begin{equation}
\|u\|_{C^{\alpha}(\tilde{U}\times(t_1,t_2))}<\infty.
\label{ugood}
\end{equation}
\newline
{\bf (B4)} For all $\phi\in C^1(U\times [0,\Lambda)\,;\,\R^+)$ with $\phi(\cdot,t)\in C^1_c(U)$ and 
$0\leq t_1<t_2<\Lambda$, we have
\begin{equation}
\|V_{t_2}\|(\phi(\cdot,t_2))-\|V_{t_1}\|(\phi(\cdot,t_1)) \leq \int_{t_1}^{t_2}{\mathcal B}(V_t,u(\cdot,t),\phi(\cdot,t))\, dt
+\int_{t_1}^{t_2}\int_U \frac{\partial \phi}{\partial t}(\cdot,t)\, d\|V_t\|dt.
\label{main}
\end{equation}
\begin{rem} As is stated in the previous section, (B4) is a weak integral form of the motion law: velocity $=$ mean curvature
$+$ $u$. In particular, if $u=0$, it is Brakke's MCF in an integral form. If there exists $\tilde{U}\subset\subset U$
such that ${\rm spt}\,\|V_t\|\subset \tilde{U}$ for all $t\in [0,\Lambda)$, then we do not need to assume (B2). In this
case, (B2) is satisfied automatically.  This can be proved easily: choose $\phi\in C^1_c(U;\R^+)$ 
with $\phi\equiv 1$ on $\tilde{U}$ and use \eqref{main} and the H\"{o}lder inequality 
to show that $\frac{d}{dt}\|V_t\|(\tilde{U})\leq \|u\|^2_{0}
\|V_t\|(\tilde{U})$, which gives a uniform bound \eqref{masft}. If we work under periodic boundary conditions (i.e., 
$U=\T^n$, for example, where $\T^n$ is the $n$-dimensional torus), we do not need (B2) by the same reason. 
\end{rem}
\subsection{Partial regularity}
\begin{define}
A point $x\in U\cap {\rm spt}\, \|V_t\|$ is said to be a $C^{2,\alpha}$ regular point if there exists 
some open neighborhood $O$ in $\R^{n}$ containing $x$
such that
$O \cap {\rm spt}\, \|V_s\| $ is an embedded 
$k$-dimensional manifold represented as the graph of $f(\cdot,s): B_R^k\rightarrow
O$ for $s\in (t-R^2,t+R^2)$ for some $R>0$ and with
\begin{equation*}
\|f\|_{0}+\|\nabla f\|_{0}+\|\nabla^2 f\|_{\alpha}+\|\p f/\p s\|_{\alpha}<\infty.
\end{equation*}
\end{define}
\begin{thm}
Under  the assumptions (B1)-(B4), for a.e$.$ $t\in (0,\Lambda)$, there exists a 
(possibly empty) closed set $G_t\subset{\rm spt}\, \|V_t\|$ with $\H^k(G_t)=0$ 
such that ${\rm spt}\, \|V_t\|\setminus G_t$ is a set of $C^{2,\alpha}$ regular points.
Moreover, we have the motion law in the classical sense, namely, the normal 
velocity vector is equal to the sum of the mean curvature vector and $u^{\perp}$ 
at each $C^{2,\alpha}$ regular point. 
\label{mainreg1}
\end{thm}
\begin{rem}
For $u=0$, Theorem \ref{mainreg1} combined with the standard linear regularity
theory proves that the above $f$ is $C^{\infty}$ on the set of $C^{2,\alpha}$ regular points.
This proves `almost everywhere regularity' of unit density Brakke's MCF. 
\end{rem}
\label{pr}
\subsection{Local regularity theorem}
To describe the local regularity theorem, we need the following (cf. \cite[Def. 5.1]{Kasai-Tonegawa})
\begin{define}
Fix $\phi\in C^{\infty}([0,\infty))$ such that $0\leq \phi\leq 1$, 
\begin{equation*}
\phi(x)\left\{\begin{array}{ll}=1 & \mbox{for }0\leq x\leq (2/3)^{1/k}, \\ >0 & \mbox{for }0\leq x<(5/6)^{1/k}, \\
=0 & \mbox{for }x\geq (5/6)^{1/k}.
\end{array}
\right. 
\end{equation*}
For $0<R<\infty$, $x\in \R^n$ and $T\in {\bf G}(n,k)$ define
\begin{equation}
\phi_{T,R}(x):=\phi(R^{-1}|T(x)|), \hspace{.2cm}
{\bf c}:=\int_{T}
\phi^2_{T,1}\, d\H^k\big(= R^{-k}\int_T\phi^2_{T,R}\, d\H^k\mbox{ for $\forall R>0$}\big).
\label{deff}
\end{equation}
\end{define}
With this we have the following
\begin{thm} Corresponding to $k,\, n$, $1\leq E_0<\infty$, $0<\nu<1$, $0<\alpha<1$, there
exist $0<\varepsilon_0<1$, $0<\sigma_0\leq 1/2$, $2<\Lambda_0<\infty$ and $1<c_{0}<\infty$
with the following property. For $T\in {\bf G}(n,k)$, $0<R<\infty$, $U=C(T,3R)$ and $(0,\Lambda)$
replaced by $(-\Lambda_0 R^2,\Lambda_0 R^2)$, 
suppose that $\{V_t\}_{-\Lambda_0 R^2\leq t\leq \Lambda_0 R^2}$ and $\{u(\cdot,t)\}_{-\Lambda_0 R^2
\leq t\leq \Lambda_0 R^2}$
satisfy (B1)-(B4). Suppose 
\begin{equation}
\sup_{-\Lambda_0 R^2\leq t\leq \Lambda_0 R^2}R^{-k} \|V_t\|(C(T,3R))\leq E_0,
\label{acon0}
\end{equation}
\begin{equation}
\mu:=\left( R^{-(k+4)}\int_{-\Lambda_0 R^2}^{\Lambda_0 R^2}\int_{C(T,3R)}|T^{\perp}(x)|^2\, d\|V_t\|dt\right)^{\frac12}
<\varepsilon_0,
\label{acon1}
\end{equation}
\begin{equation}
\|u\|_{\alpha}:=R\|u\|_{0}+R^{1+\alpha}[ u ]_{\alpha}<\varepsilon_0, 
\label{acon2}
\end{equation}
\begin{equation}
(-\Lambda_0+3/2)R^2 \leq \exists t_1\leq (-\Lambda_0+2) 
R^2\hspace{.2cm}:\hspace{.2cm}R^{-k}\|V_{t_1}\|(\phi^2_{T,R})<(2-\nu){\bf c},
\label{acon3}
\end{equation}
\begin{equation}
(\Lambda_0-2)R^2\leq \exists t_2\leq (\Lambda_0-3/2)R^2\hspace{.2cm}:\hspace{.2cm} R^{-k}
\|V_{t_2}\|(\phi^2_{T,R})>\nu {\bf c}.
\label{acon4}
\end{equation}
Denote $D:=(T\cap B_{\sigma_0 R})\times(-R^2 /4, R^2/4)$. Then there are $f:D\rightarrow
T^{\perp}$ and $F:D\rightarrow \R^n$ such that $T(F(y,t))=y$ and $T^{\perp}(F(y,t))=f(y,t)$
for all $(y,t)\in D$,
\begin{equation}
{\rm spt}\, \|V_t\|\cap C(T,\sigma_0 R)={\rm image}\, F(\cdot,t)\hspace{.3cm}\forall t\in (-R^2/4,R^2/4),
\label{acon5}
\end{equation}
\begin{equation}
\mbox{$f$ is twice differentiable w.r.t. $x$ and differentiable w.r.t. $(x,t)$ on $D$,}
\label{acon6}
\end{equation}
\begin{equation}
R\Big\|R^{-2}f,\,R^{-1}\nabla f,\,\nabla^2 f,\,\frac{\p f}{\p t}\Big\|_{0}+R^{1+\alpha}\Big[\nabla^2 f,\,\frac{\p f}{\p t}\Big]_{\alpha}
\leq c_0\max\{
\mu,\|u\|_{\alpha}\}.
\label{acon7}
\end{equation}
Moreover the motion law (normal 
velocity $=$ mean curvature vector $+\, u^{\perp}$) is satisfied on ${\rm image}\, F$. 
\label{mainreg2}
\end{thm}
\eqref{acon1} requires smallness of deviation from $k$-dimensional 
plane in a weak measure-theoretic sense. \eqref{acon3} excludes the possibility that
there may be two or more almost parallel $k$-dimensional planes which may not 
move for the whole time. Obviously, for such case, we cannot hope to represent the
graph as a univalent function. The idea of having possibly large $\Lambda_0$ is that,
if we have a mass strictly less than that of 2 sheets of $k$-dimensional planes near the
beginning, we will have a nice univalent representation of graph after sufficiently 
long time. Asking a certain mass lower bound \eqref{acon4} is also natural since $V_t=0$ for all time 
would satisfy (B1)-(B4) as well as \eqref{acon0}-\eqref{acon3}. Since one can always set $V_t=0$
after any instance and still obtain a solution satisfying \eqref{main}, we need to impose \eqref{acon4} 
towards the end of the time interval. 
\section{$L^2$-$L^{\infty}$ estimate}
In this section we first define function $Q_g$, which is a (square of)
distance function from a graph of solution of the heat equation, 
roughly speaking. We then prove that the 
$L^2$ norm of $Q_g$ controls such distance function in sup-norm in Proposition \ref{height}, 
which is analogous to $L^2$-$L^{\infty}$ estimate of \cite[Prop. 6.4]{Kasai-Tonegawa}.
The Dirichlet energy of the distance is also similarly controlled. 
Throughout this section let $T\in {\bf G}(n,k)$ be the projection matrix corresponding to ${\mathbb R}^k\times\{0\}$. 

Suppose that we are given a function $g=(g_{k+1},\cdots,g_n)$ defined on ${\mathbb R}^n\times{\mathbb R}$ with
the following conditions.
For each $l=k+1,\cdots,n$, 
\begin{equation}
g_l(x_1,\cdots,x_n,t)=a_l+b_l  t+\sum_{i=1}^k a_{li} x_i+\frac12\sum_{i,j=1}^k a_{lij}x_i x_j
\label{propfl1}
\end{equation}
for some $a_l,\, b_l,\, a_{li},\, a_{lij}\in {\mathbb R}$ with $a_{lij}=a_{lji}$ for all $1\leq i,j\leq k$.
Note that $g_l$ depends only on $t$ and $x_1,\cdots, x_k$ and we often consider $g_l$ as a function defined
on ${\mathbb R}^k\times {\mathbb R}$.  We additionally assume that 
\begin{equation}
b_l=\sum_{i=1}^k a_{lii}. 
\label{propfl2}
\end{equation}
Equivalently, each component function $g_l$ of $g$ satisfies the heat equation 
$\frac{\partial g_l}{\partial t}=\Delta g_l$. 
We next define
\begin{define}
If $g=(g_{k+1},\cdots,g_n)$ satisfies \eqref{propfl1} and \eqref{propfl2}, then we write
\[g\in \F.\]
For $g\in \F$, we define a function $Q_g$ defined on ${\mathbb R}^n\times{\mathbb R}$ by
\begin{equation}
Q_g(x,t):=\frac12\sum_{l=k+1}^{n}(x_{l}-g_l(x,t))^2.
\label{defQ}
\end{equation}
\end{define}
Note that $(2Q_g)^{1/2}$ is the vertical distance of the 
point $x$ from the graph of $g$. The expectation is that the MCF
should be closely approximated by the solution of the heat equation. We next need the following 
technical lemma.
\begin{lemma}
There exists $1<c_1(n,k)<\infty$ 
with the following property. Suppose 
a function $f=(f_{k+1},\cdots,f_{n})
:B_1^k\times(-1,1)\rightarrow {\mathbb R}^{n-k}$ with continuous $\nabla f$ is given. Define $M_t:={\rm graph}\, f(\cdot,t)\subset {\mathbb R}^n$. 
We assume that 
\begin{equation}
\sup_{B_1^k\times(-1,1)} |\nabla f |\leq 1.
\label{gracon}
\end{equation}
Suppose $g=(g_{k+1},\cdots,g_n) \in \F$
is given and
for each $l=k+1,\cdots,n$ define 
${\bf g}_l:B_1^k\times(-1,1)\rightarrow {\mathbb R}^n$ by 
\begin{equation}
{\bf g}_l:=\big(\frac{\partial g_l}{\partial x_1},\cdots,
\frac{\partial g_l}{\partial x_k},0,\cdots,-1,\cdots, 0\big),
\label{defF}
\end{equation} 
where $-1$ is in the $l$-th component of ${\bf g}_l$.
For each $(x,t)\in B_1^k\times(-1,1)$ let $S=S(x,t)\in {\bf G}(n,k)$ be the tangent space $T_{(x,f(x,t))} M_t$. Then we have
\begin{equation}
\frac{\partial Q_g}{\partial t}-S\cdot \nabla^2 Q_g\leq c_1 Q_g^{1/2} |\nabla f|^2 |\nabla^2 g| - \frac{1}{4k} |\nabla f-\nabla g|^2-\frac12\sum_{l=k+1}^n |S({\bf g}_l)|^2
\label{Q2d}
\end{equation}
and
\begin{equation}
\frac{\partial g_l}{\partial t}-S\cdot\nabla^2 g_l\leq c_1 |\nabla f|^2 |\nabla ^2 g_l|.
\label{Q2e}
\end{equation}
Note that $\frac{\p Q_g}{\p t}$, $\nabla^2 Q_g$ and $Q_g^{1/2}$ are evaluated at $(x,f(x,t))
\in M_t$ in \eqref{Q2d}.
\label{propQ}
\end{lemma}
{\it Proof}.  One checks that 
\begin{equation}
\nabla^2 Q_g=\sum_{l=k+1}^n ({\bf g}_l\otimes {\bf g}_l- (x_l-g_l)\nabla^2 g_l)
\label{Q2}
\end{equation}
where $\nabla^2 g_l$ is the $n\times n$ matrix with non-zero components only in the upper-left
$k\times k$ sub-matrix. 
Due to 
\eqref{propfl2} and \eqref{Q2}, we have
\begin{equation}
\frac{\partial Q_g}{\partial t}-S\cdot\nabla^2 Q_g=\sum_{l=k+1}^n (-S\cdot({\bf g}_l\otimes {\bf g}_l)+(x_l-g_l)(T-S)\cdot \nabla^2 g_l).
\label{Q2a}
\end{equation}
We estimate each term of the right-hand side of \eqref{Q2a}. For the first term, since $S\in {\bf G}(n,k)$, we have
$S\cdot({\bf g}_l\otimes {\bf g}_l)=|S({\bf g}_l)|^2$. Fix any $l=k+1,\cdots,n$ and $j=1,\cdots, k$. 
Since $S$ is the tangent space of ${\rm graph}\,(f_{k+1},\cdots,f_n)$, 
note that $S$ contains ${\bf f}_{j}:=(0,\cdots, 1,\cdots,0,\frac{\partial f_{k+1}}{\partial x_{j}},\cdots,\frac{\partial f_n}{\partial x_{j}})$,
where $1$ is in the $j$-th component of ${\bf f}_j$. Thus we may conclude that
\begin{equation}
S\cdot({\bf g}_l\otimes {\bf g}_l)=|S({\bf g}_l)|^2\geq |{\bf g}_l\cdot {\bf f}_j|^2/|{\bf f}_j|^2 =\Big|\frac{\partial g_l}{\partial x_j}-\frac{\partial f_l}{\partial x_j}\Big|^2
\Big/\Big(1+\Big|\frac{\partial f}{\partial x_j}\Big|^2\Big).
\label{Q2b}
\end{equation}
By \eqref{gracon} and summing over $j$ and $l$, we obtain from \eqref{Q2b}
\begin{equation}
\sum_{l=k+1}^n S\cdot({\bf g}_l\otimes {\bf g}_l)\geq \frac{1}{2k} |\nabla f-\nabla g|^2.
\label{Q2c}
\end{equation}
In particular, from \eqref{Q2c}, we obtain
\begin{equation}
\sum_{l=k+1}^n S\cdot({\bf g}_l\otimes {\bf g}_l)\geq\frac{1}{4k}|\nabla f-\nabla g|^2+ \frac12 \sum_{l=k+1}^n |S({\bf g}_l)|^2.
\label{Q2cs}
\end{equation}
For the second term of \eqref{Q2a}, we need to know the expression of $T-S$.
The $k$-dimensional space corresponding to $S$ is spanned by ${\bf f}_1,\cdots,{\bf f}_k$. 
Consider the Gram-Schmidt orthonormalization $\tilde{{\bf f}}_1,\cdots, \tilde{{\bf f}}_k$ of 
${\bf f}_1,\cdots,{\bf f}_k$, namely, $\tilde{{\bf f}}_1={\bf f}_{1}/|{\bf f}_1|$, $\hat{{\bf f}}_2={\bf f}_2-({\bf f}_2\cdot \tilde{{\bf f}}_1)
\tilde{{\bf f}}_1$, $\tilde{{\bf f}}_2=\hat{{\bf f}}_2/|\hat{{\bf f}}_2|$, $\cdots$, $\hat{{\bf f}}_k={\bf f}_k-\sum_{j=1}^{k-1}
({\bf f}_k\cdot \tilde{{\bf f}}_j)\tilde{{\bf f}}_j$, $\tilde{{\bf f}}_k=\hat{{\bf f}}_k/|\hat{{\bf f}}_k|$. Then $S=\sum_{j=1}^k \tilde{{\bf f}}_j\otimes\tilde{{\bf f}}_j$. 
It is not difficult to check that each entry of the upper-left $k\times k$ sub-matrix of 
$T-S$ is bounded by some constant times $|\nabla f|^2$, where the constant depends 
only on $k$ and $n$. The reason is as follows. 
The first $k$ components of $\hat{{\bf f}}_j$ are $(O(|\nabla f|^2),\cdots,O(|\nabla f|^2),1,0,\cdots, 0)$, 
where $1$ is in the $j$-th component. The last $n-k$ components of $\hat{{\bf f}}_j$ are $O(|\nabla f|)$. The division
by $1/\sqrt{1+O(|\nabla f|^2)}$ for normalization does not change the order of magnitude except that the $j$-th 
component turns $1+O(|\nabla f|^2)$. One sees that the next vector $\tilde{{\bf f}}_{j+1}$ has the 
same property. Thus for each $j=1,\cdots,k$, $\tilde{{\bf f}}_j\otimes \tilde{{\bf f}}_j$ has $O(|\nabla f|^2)$ components
for the upper-left $k\times k$ sub-matrix except for the $j-j$ component, which is $1+O(|\nabla f|^2)$. 
Since $T$ has $1$ in the diagonal components for the upper-left $k\times k$ sub-matrix, we 
have the above stated property. Note that we only need to consider such entries since $\nabla^2 g_l$
has non-zero entries only there. Thus with \eqref{Q2a} and 
\eqref{Q2cs}, we obtain \eqref{Q2d}. The derivation for \eqref{Q2e} is similar, which only requires the estimate
for $(T-S)\cdot\nabla^2 g_l$. 
\hfill{$\Box$}
\begin{prop}
There exists $c_2=c_2(n,k)$ with the following property. 
Suppose that $\{V_t\}_{-1<t<1}$ and $\{u(\cdot, t)\}_{-1<t<1}$, where $V_t=|M_t|$ with 
$M_t={\rm graph}\, f(\cdot,t)$, satisfy (B1) and (B4) on $C(T,1)\times (-1,1)$.
Let $g\in \F$ be given with $Q_g$ as in \eqref{defQ}. 
In addition, assume \eqref{gracon},
\begin{equation}
\sup_{B_1^k\times(-1,1)} |\nabla g|\leq 1
\label{amoQ}
\end{equation}
and
\begin{equation}
\|u\|_{0}:=\sup_{C(T,1)\times (-1,1)}|u|\leq 1.
\label{amou}
\end{equation}
Then we have 
\begin{equation}
\begin{split}
&\sup_{B_{1/2}^k\times(-3/4,1)}|f-g|^2
+\int_{-3/4}^1\int_{B_{1/2}^k}|\nabla f-\nabla g|^2\, d{\mathcal H}^k dt
\\ & \leq c_2 \big(\int_{-1}^1\int_{C(T,1)} Q_g\, d\|V_t\|dt+  \|u\|_{0}^2+\|\nabla f\|_{0}^4 \|\nabla^2 g\|^2_{0}\big).
\end{split}
\label{con1}
\end{equation}
\label{height}
\end{prop}
{\it Proof}. In the proof let $\eta\in C^{\infty}(B_1^k\times (-1,1))$ be a non-negative function with $\eta\equiv 1$ on 
$B_{3/4}^k\times(-7/8,1)$, $\eta\equiv 0$ on $B_1^k\times(-1,1)\setminus B_{7/8}^k\times(-15/16,1)$, $0\leq \eta\leq 1$
and $|\nabla\eta|,\, |\nabla^2 \eta|,\, |\frac{\partial\eta}{\partial t}|\leq c(k)$. We then re-define $\eta(x,t):=\eta(T(x),t)$ for
$(x,t)\in C(T,1)\times(-1,1)$. 
For $(y,s)\in C(T,1/2)\times (-3/4,\infty)$, we use $\phi(x,t)=Q_g(x,t)\rho_{(y,s)}(x,t)\eta(x,t)$ in \eqref{main}, over the 
time interval $t_1=-1$ and $-1<t_2<\min\{s,1\}$. We then obtain (writing $\rho_{(y,s)}(x,t)$ as $\rho$ and $Q_g(x,t)$ 
as $Q$)
\begin{equation}
\left.\int_{C(T,1)} Q\rho\eta\, d\|V_t\|\right|_{t=t_2}
\leq \int_{-1}^{t_2}\int_{C(T,1)}\{-h\rho\eta Q+\nabla(\rho\eta Q)\}\cdot(h+u^{\perp})+\frac{\partial}{\partial t}(Q\rho\eta)
\, d\|V_t\|dt
\label{moQ1}
\end{equation}
since $\eta=0$ for $t=-1$. For a.e. $t\in (-1,t_2)$, we may compute the integrand of the right-hand side of \eqref{moQ1} 
as follows. Here we use the perpendicularity of mean curvature \eqref{perpH} in deriving $\nabla\rho\cdot h=(\nabla\rho)^{\perp}
\cdot h$.
\begin{equation}
\begin{split}
-|h|^2 \rho\eta Q&+(\nabla\rho\cdot h)\eta Q+\rho\nabla(\eta Q)\cdot h+(-h\rho\eta Q+\eta Q \nabla\rho)\cdot u^{\perp}
+\rho\nabla(\eta Q)\cdot u^{\perp}+\frac{\partial}{\partial t}(\rho\eta Q)\\
&\leq -\rho\big|h-\frac{(\nabla\rho)^{\perp}}{\rho}\big|^2\eta Q-(\nabla\rho\cdot h)\eta Q+\frac{|(\nabla\rho)^{\perp}|^2}{\rho}
\eta Q +\rho\nabla(\eta Q)\cdot h \\
&+\rho\big|h-\frac{(\nabla\rho)^{\perp}}{\rho}\big|^2\eta Q+\rho\eta Q|u|^2+\rho
\nabla(\eta Q)\cdot u^{\perp}+\frac{\partial}{\partial t}(\rho\eta Q).
\end{split}
\label{moQ2}
\end{equation}
Thus we have 
\begin{equation}
\begin{split}
\left.\int_{C(T,1)}Q\rho\eta\, d\|V_t\|\right|_{t=t_2}&\leq \int_{-1}^{t_2}\int_{C(T,1)}-(\nabla\rho\cdot h)\eta Q+\rho\nabla(\eta Q)\cdot h
+\frac{|(\nabla\rho)^{\perp}|^2}{\rho}\eta Q\\
&+\rho\eta Q|u|^2 +\rho\nabla(\eta Q)\cdot u^{\perp}+\frac{\partial}{\partial t}(\eta\rho Q)\, d\|V_t\|dt.
\end{split}
\label{moQ3}
\end{equation}
By \eqref{first}, the first two terms of the right-hand side of \eqref{moQ3} is
\begin{equation}
\begin{split}
\int_{-1}^{t_2}\int_{G_k(C(T,1))}&\nabla(\eta Q\nabla \rho)\cdot S-\nabla\{\rho\nabla(\eta Q)\}\cdot S\, dV_t(\cdot,S)dt\\
&=\int_{-1}^{t_2}\int_{G_k(C(T,1))}(\nabla^2 \rho\cdot S)\eta Q-\rho\nabla^2(\eta Q)\cdot S\, dV_t(\cdot, S)dt.
\end{split}\label{moQ4}
\end{equation}
Using 
\begin{equation}
\nabla^2\rho\cdot S+\frac{|(\nabla\rho)^{\perp}|^2}{\rho}+\frac{\partial\rho}{\partial t}=0,
\label{moQ5sub}
\end{equation}
we obtain from
\eqref{moQ3} and \eqref{moQ4}
\begin{equation}
\begin{split}
&\left.\int_{C(T,1)}Q\rho\eta\, d\|V_t\|\right|_{t=t_2}\leq \int_{-1}^{t_2}\int_{G_k(C(T,1))}\{-\nabla^2 (\eta Q)\cdot S+\frac{\partial}{\partial t}(Q\eta)
\\&+\eta Q |u|^2+Q\nabla\eta\cdot u^{\perp}+\eta\nabla Q\cdot u^{\perp}\}\rho\, dV(\cdot,S)dt=: I_1+I_2+I_3+I_4+I_5.
\end{split}
\label{moQ5}
\end{equation}
{\it Estimate of $I_1+I_2$}. The integrand of $I_1$ is
\begin{equation}
\{-Q S\cdot \nabla^2 \eta-2 (\nabla \eta\otimes \nabla Q)\cdot S-\eta\nabla^2 Q\cdot S\}\rho.
\label{moQ6}
\end{equation}
Note that, with the notation of \eqref{defF}, we have
\begin{equation}
(\nabla \eta\otimes \nabla Q)\cdot S=\nabla\eta\cdot S(\nabla Q)=-\nabla\eta\cdot\sum_{l=k+1}^n
(x_l-g_l)S({\bf g}_l).
\label{moQ7}
\end{equation}
Thus, we obtain from \eqref{moQ7} and the Cauchy-Schwarz inequality that
\begin{equation}
-2(\nabla\eta\otimes\nabla Q)\cdot S\leq 2\sqrt{2}|\nabla\eta| Q^{\frac12}\big(\sum_{l=k+1}^n|S({\bf g}_l)|^2\big)^{\frac12}
\leq \frac12 \sum_{l=k+1}^n |S({\bf g}_l)|^2\eta+4Q\frac{|\nabla\eta|^2}{\eta}.
\label{moQ8}
\end{equation}
By \eqref{moQ6}, \eqref{moQ8} and Lemma \ref{propQ}, we obtain
\begin{equation}
\begin{split}
I_1+I_2\leq \int_{-1}^{t_2}&\int_{C(T,1)}Q\rho\big(|\nabla^2\eta|+4\frac{|\nabla\eta|^2}{\eta}+\big|\frac{\partial \eta}{\partial t}\big|\big) \\&+
c_1 Q^{\frac12}\rho|\nabla f|^2 |\nabla^2 g|\eta-\frac{1}{4k}|\nabla f-\nabla g|^2\eta\rho\, d\|V_t\|dt.
\end{split}
\label{moQ9}
\end{equation}
For $\rho=\rho_{(y,s)}(x,t)$ with $(y,s)\in C(T,1/2)\times(-3/4,\infty)$ and $(x,t)\notin \{\eta=1\} $
we have $|x-y|\geq 1/4$ or $s-t\geq 1/8$. Thus we have
\begin{equation}
\rho(|\nabla^2\eta|+4\frac{|\nabla\eta|^2}{\eta}+|\frac{\partial \eta}{\partial t}|)\leq c(k)
\label{moQ10}
\end{equation}
for a suitable constant depending only on $k$. Then by the Cauchy-Schwarz inequality, \eqref{moQ9} and \eqref{moQ10} give
\begin{equation}
I_1+I_2\leq \int_{-1}^{t_2}\int_{C(T,1)}c(k)Q+Q\eta\rho+c_1^2|\nabla f|^4|\nabla^2 g|^2 \eta\rho
-\frac{1}{4k}|\nabla f-\nabla g|^2\rho\eta\, d\|V_t\|dt.
\label{moQ11}
\end{equation}
Since $\int_{C(T,1)}\rho\eta\, d\|V_t\|\leq c(k)$ by \eqref{gracon}, we obtain from \eqref{moQ11}
\begin{equation}
I_1+I_2\leq \int_{-1}^{t_2}\int_{C(T,1)}c(k)Q+Q\eta\rho-\frac{1}{4k}|\nabla f-\nabla g|^2\rho\eta\, d\|V_t\|dt
+2c_1^2 c(k)\|\nabla f\|_{0}^4 \|\nabla^2 g\|_{0}^2.
\label{moQ12}
\end{equation}
{\it Estimate of $I_3$}. We have by \eqref{amou}
\begin{equation}
\int_{-1}^{t_2}\int_{C(T,1)}\eta Q |u|^2\rho\, d\|V_t\|dt\leq \int_{-1}^{t_2}\int_{C(T,1)}\eta Q\rho\, d\|V_t\|dt.
\label{moQ13}
\end{equation}
{\it Estimate of $I_4$}. By \eqref{amou} and since $\rho\leq c(k)$ on the support of $|\nabla\eta|$, we have for any $t_2\in (-1,1)$
\begin{equation}
\int_{-1}^{t_2}\int_{C(T,1)}Q\nabla\eta\cdot u^{\perp}\rho\, d\|V_t\|dt\leq c(k)\int_{-1}^1\int_{C(T,1)}Q\, d\|V_t\|dt.
\label{moQ14}
\end{equation}
{\it Estimate of $I_5$}. By \eqref{amoQ}, one can check that $|\nabla Q|\leq \sqrt{Q} c(k)$, thus
\begin{equation}
\int_{-1}^{t_2}\int_{C(T,1)} \eta \nabla Q\cdot u^{\perp}\rho\, d\|V_t\|dt\leq 
c(k)\|u\|_{0}^2 +\int_{-1}^{t_2}\int_{C(T,1)}\eta Q\rho\, d\|V_t\|dt.
\label{moQ15}
\end{equation}
Thus, with a suitable $c=c(k)$, we have from \eqref{moQ5}, \eqref{moQ12}, \eqref{moQ13}-\eqref{moQ15}
\begin{equation}
\begin{split}
&\left.\int_{C(T,1)} Q\rho\eta\, d\|V_t\|\right|_{t=t_2} \leq c\big(\int_{-1}^1\int_{C(T,1)}Q\, d\|V_t\|dt
+\int_{-1}^{t_2}\int_{C(T,1)}Q\eta\rho\, d\|V_t\|dt \\
& + \|u\|_{0}^2+ \|\nabla f\|_{0}^4\|\nabla^2 g\|_{0}^2\big)
-\frac{1}{4k}\int_{-1}^{t_2}\int_{C(T,1)}|\nabla f-\nabla g|^2 \eta\rho\, d\|V_t\|dt
\end{split}
\label{moQ16}
\end{equation}
for all $t_2\in (-1,\min\{s,1\})$. Note that the differential inequality $F'\leq c(F+\tilde{c})$ with $F(-1)=0$ implies 
$F(t)\leq \tilde{c}(e^{c(t+1)}-1)$ and thus $F'(t)\leq c\tilde{c}e^{c(t+1)}$.  Thus by dropping the last term of \eqref{moQ16}, we obtain 
\begin{equation}
\left.\int_{C(T,1)}Q\eta\rho\, d\|V_t\|\right|_{t=t_2}\leq c e^{2c} \big(\int_{-1}^1 \int_{C(T,1)}Q\, d\|V_t\|dt+\|u\|_{0}^2
+\|\nabla f\|_{0}^4 \|\nabla^2 g\|_{0}^2\big).
\label{moQ17}
\end{equation}
For any $t_2\in (-3/4,1)$ and $y\in C(T,1/2)\cap {\rm spt}\, \|V_{t_2}\|$, and arbitrarily small $\varepsilon>0$, 
we use $\rho=\rho_{(y,t_2+\varepsilon)}(x,t)$ in the above computation. Since $\rho_{(y,t_2+\varepsilon)}(\cdot,
t_2)\|V_{t_2}\|\rightharpoonup \delta_{y}$ as $\varepsilon\rightarrow 0$, where $\delta_{y}$ is the delta
function at $y$, and since $\eta=1$ on $C(T,1/2)\times(-3/4,1)$, we conclude that the first term of the left-hand side
of \eqref{con1} is bounded by the right-hand side of \eqref{moQ17}. For the second term of the right-hand side of
\eqref{con1}, we use \eqref{moQ16} with $\rho\equiv 1$. Note that the only property we used for the above 
computation is \eqref{moQ5sub}. Since
\begin{equation}
\int_{-3/4}^1\int_{B_{1/2}^k}|\nabla f-\nabla g|^2\, d{\mathcal H}^k dt\leq 
\int_{-3/4}^1\int_{C(T,1)}|\nabla f-\nabla g|^2\, d\|V_t\|dt,
\label{moQ18}
\end{equation}
we prove \eqref{con1}.
\hfill{$\Box$}
\section{A decay estimate by blow-up argument}
The main result of this section is the following Proposition
\ref{blow} which shows that one can find a better approximation in $\F$
in a smaller scale if the relevant quantities \eqref{b1}-\eqref{b4} are sufficiently small.
It is similar to \cite[Proposition 8.1]{Kasai-Tonegawa}, except that the decay we obtain
here is $\theta^{1+\a}$ instead of $\theta^{\varsigma}$. 
\begin{prop}
Corresponding to $n$, $k$ and $0<\alpha<1$ there exist $0<\varepsilon_1<1$,
$0<\theta<1/4$, $1<c_3<\infty$ with the following property. For $0<R<\infty$, suppose
$\{V_t\}_{-R^2<t<R^2}$ and $\{u(\cdot,t)\}_{-R^2<t<R^2}$, where $V_t=|M_t|$ with
$M_t={\rm graph}\, f(\cdot,t)$, satisfy (B1) and (B4) on $C(T,R)\times
(-R^2,R^2)$ and $g\in \F$ is given. Denote $\|\cdot\|_{0}:=
\sup_{B_R^k\times(-R^2,R^2)} |\cdot |$. Assume that
\begin{equation}
\|\nabla f\|_{0}\leq \varepsilon_1,
\label{b1}
\end{equation}
\begin{equation}
\|\nabla g\|_{0}+
R\|\nabla^2 g\|_{0}+R\|\frac{\partial g}{\partial t}\|_{0}\leq \varepsilon_1,
\label{b2}
\end{equation}
\begin{equation}
u(0,0)=0,
\label{b3}
\end{equation}
\begin{equation}
\mu:=\Big(R^{-(k+4)}\int_{-R^2}^{R^2}\int_{C(T,R)} Q_g\, d\|V_t\|dt\Big)^{\frac12}\leq \varepsilon_1.
\label{b4}
\end{equation}
Then there exists $\hat{g}\in \F$ with
\begin{equation}
R^{-1}\|g-\hat{g}\|_{0}+\|\nabla g-\nabla \hat{g}\|_{0}+
R\|\nabla^2 g-\nabla^2\hat{g}\|_{0}+R\|\frac{\partial g}{\partial t}-\frac{\partial\hat{g}}{\partial t}\|_{0}\leq c_3\mu
\label{b5}
\end{equation}
and 
\begin{equation}
\Big((\theta R)^{-(k+4)}  \int_{-\theta^2 R^2}^{\theta^2 R^2}\int_{C(T,\theta R)} Q_{\hat{g}}\, d\|V_t\|dt\Big)^{\frac12} 
\leq \theta^{1+\alpha} \max\{\mu,c_3 R^{1+\alpha}[ u ]_{\alpha},c_3  \|\nabla f\|_{0}^2\}.
\label{b6}
\end{equation}
\label{blow}
\end{prop}
{\it Proof}. After a change of variables, we may assume $R=1$. Note that the statement is written in a scale invariant manner. 
If the claim were false, then for each $m\in {\mathbb N}$ there exist $\{V_t^{(m)}\}_{-1<t<1}$ (represented by $f^{(m)}$), 
$\{u^{(m)}(\cdot,t)\}_{-1<t<1}$ satisfying \eqref{b3}, (B1) and (B4) on $C(T,1)\times (-1,1)$ and $g^{(m)}\in \F$ such that 
\begin{equation}
\|\nabla f^{(m)}\|_{0}\leq \frac{1}{m},
\label{b7}
\end{equation}
\begin{equation}
\|\nabla g^{(m)}\|_{0}+
\|\nabla^2 g^{(m)}\|_{0}+\|\frac{\partial g^{(m)}}{\partial t}\|_{0}\leq \frac{1}{m},
\label{b8}
\end{equation}
\begin{equation}
\mu^{(m)}:=\left(\int_{-1}^1\int_{C(T,1)} Q_{g^{(m)}}\, d\|V_t^{(m)}\|dt\right)^{\frac12}\leq \frac{1}{m},
\label{b9}
\end{equation}
but for any $g\in \F$ with
\begin{equation}
\|g-g^{(m)}\|_{0}+\|\nabla g-\nabla g^{(m)}\|_{0}+
\|\nabla^2 g-\nabla^2 g^{(m)}\|_{0}+\|\frac{\partial g}{\partial t}-\frac{\partial g^{(m)}}{\partial t}\|_{0}\leq m\mu^{(m)},
\label{b10}
\end{equation}
we have
\begin{equation}
\left(\theta^{-(k+4)}\int_{-\theta^2}^{\theta^2}\int_{C(T,\theta)} Q_g\, d\|V_t^{(m)}\|dt\right)^{\frac12}>\theta^{1+\alpha}
\max\{\mu^{(m)},m[ u^{(m)}]_{\alpha},m\|\nabla f^{(m)}\|_{0}^2 \}.
\label{b11}
\end{equation}
Here, $0<\theta<1/4$ will be chosen depending only on $k$, $n$ and $\alpha$. 
By using $g=g^{(m)}$ (which satisfies \eqref{b10} trivially), we obtain from \eqref{b11}
\begin{equation}
\max\{[u^{(m)}]_{\alpha}, \|\nabla f^{(m)}\|_{0}^2\}\leq \frac{\theta^{-\frac{k+4}{2}-1-\alpha}\mu^{(m)}}{m}.
\label{b12}
\end{equation}
Thus \eqref{b12} shows
\begin{equation}
\lim_{m\rightarrow\infty}(\mu^{(m)})^{-1}[u^{(m)}]_{\alpha}=\lim_{m\rightarrow\infty}(\mu^{(m)})^{-1}\|\nabla f^{(m)}\|_{0}^2
=0.
\label{b13}
\end{equation}
Next we use Proposition \ref{height} for $f=f^{(m)}$, $g=g^{(m)}$, $u=u^{(m)}$.  
The required conditions \eqref{gracon}, \eqref{amoQ} and 
\eqref{amou} follow from \eqref{b7}, \eqref{b8}, \eqref{b3} and \eqref{b13}. Then we have
\begin{equation}
\sup_{B_{1/2}^k\times(-3/4,1)} |f^{(m)}-g^{(m)}|^2
+\int_{-3/4}^1\int_{B_{1/2}^k}|\nabla f^{(m)}-\nabla g^{(m)}|^2\, d{\mathcal H}^k dt\leq 2 c_2 (\mu^{(m)})^2
\label{b14}
\end{equation}
for all sufficiently large $m$ due to \eqref{con1}, \eqref{b3}, \eqref{b8} and \eqref{b13}. We next define a sequence of renormalized functions
\begin{equation}
\tilde{f}^{(m)}(x,t):=(\mu^{(m)})^{-1}(f^{(m)}(x,t)-g^{(m)}(x,t))
\label{b15}
\end{equation}
with $\tilde{f}^{(m)}=(\tilde{f}^{(m)}_{k+1},\cdots,\tilde{f}^{(m)}_{n})$. 
From \eqref{b14} and \eqref{b15}, we have
\begin{equation}
\sup_{B_{1/2}^k\times(-3/4,1)} |\tilde{f}^{(m)}|^2+\int_{-3/4}^1\int_{B_{1/2}^k}|\nabla\tilde{f}^{(m)}|^2\, d{\mathcal H}^k dt\leq 2c_2
=: (c_4)^2
\label{b16} 
\end{equation}
for all sufficiently large $m$. By the standard compactness theorem, 
there exist a convergent subsequence (denoted by the same index) and a limit $\tilde{f}=(\tilde{f}_{k+1},\cdots,\tilde{f}_n)$ such that 
\begin{equation}
\tilde{f}^{(m)}\rightharpoonup \tilde{f}\hspace{.3cm}\mbox{weakly in }L^2,\hspace{1cm}
\nabla\tilde{f}^{(m)}\rightharpoonup\nabla\tilde{f}\hspace{.3cm}\mbox{weakly in }L^2
\label{b17}
\end{equation}
both on $B_{1/2}^k\times (-3/4,1)$ and
\begin{equation}
 \|\tilde{f}\|_{L^{\infty}(B_{1/2}^k\times(-3/4,1))}^2+\int_{-3/4}^1\int_{B_{1/2}^k}|\nabla\tilde{f}|^2\, d{\mathcal H}^k dt\leq (c_4)^2
\label{b17.1}
\end{equation}
In addition, due to Rellich's compactness theorem, we may choose such 
subsequence so that for a countable dense set $\{s_j\}_{j=1}^{\infty}\subset (-3/4,1)$, 
\begin{equation}
\{\tilde{f}^{(m)}(\cdot,s_j)\}_{m=1}^{\infty}\,\mbox{ is a Cauchy sequence in $L^2 (B_{1/2}^k)$
for all $j\in {\mathbb N}$. }
\label{b17.5}
\end{equation}
We next claim
\begin{lemma}
Each component function $\tilde{f}_l$ is in $C^{\infty}(B_{1/2}^k\times(-3/4,1))$ and satifies the heat equation,
\begin{equation}
\frac{\partial \tilde{f}_l}{\partial t}-\Delta \tilde{f}_l=0
\label{b18}
\end{equation}
on $B_{1/2}^k\times(-3/4,1)$.
\label{heateq}
\end{lemma}
{\it Proof of Lemma \ref{heateq}}. In the following we fix $l\in \{k+1,\cdots,n\}$. 
Let $\phi\in C_c^{\infty}(B_{1/2}^k\times(-3/4,1);\R^+)$ be arbitrary and fixed. 
For $(x,t)\in C(T,1/2)\times(-3/4,1)$ and $m\in \N$ define a function
\begin{equation}
\phi^{(m)}(x,t):=(x_l-g_l^{(m)}(T(x),t)+2c_4 \mu^{(m)})\phi(T(x),t).
\label{b19}
\end{equation}
Since $x_l=f^{(m)}_l(x,t)$ for $x\in {\rm spt}\, \|V^{(m)}_t\|$, note that we have
from \eqref{b15} and \eqref{b19}
\begin{equation}
\phi^{(m)}(x,t)=\mu^{(m)}(\tilde{f}_l^{(m)}(T(x),t)+2c_4)\phi(T(x),t)
\label{b20}
\end{equation}
for $x\in {\rm spt}\, \|V_t^{(m)}\|$. 
Thus, due to \eqref{b16}, $\phi^{(m)}$ is non-negative on ${\rm spt}\, \|V_t^{(m)}\|$ for all
$t\in (-3/4,1)$. Away from $\cup_{t\in (-3/4,1)}(C(T,1/2)\cap{\rm spt}\, \|V_t^{(m)}\|)\times\{t\}$, we may 
modify $\phi^{(m)}$ so that the modified $\phi^{(m)}$ is non-negative smooth function with
compact support in $C(T,1/2)\times(-3/4,1)$. This modification justifies the use of $\phi^{(m)}$ in
\eqref{main} but does not affect the following computations since only 
the values in some neighborhood of ${\rm spt}\, \|V^{(m)}_t\|$ matter. 
With this modification, the substitution of 
$\phi^{(m)}$ in \eqref{main} gives (denoting $h(V_t^{(m)},\cdot)$ by $h^{(m)}$) 
\begin{equation}
0\leq \int_{-3/4}^1\int_{C(T,1/2)}(-h^{(m)}\phi^{(m)}+\nabla\phi^{(m)})\cdot
(h^{(m)}+(u^{(m)})^{\perp})+\frac{\p \phi^{(m)}}{\p t}\, d\|V_t^{(m)}\|dt.
\label{b21}
\end{equation}
By the Cauchy-Schwarz inequality and dropping a negative term, \eqref{b21} gives
\begin{equation}
\begin{split}
0 \leq &\int_{-3/4}^1\int_{C(T,1/2)}|u^{(m)}|^2\phi^{(m)}+|u^{(m)}||\nabla\phi^{(m)}|+\mu^{(m)}(\tilde{f}_l^{(m)}+2c_4)
\nabla\phi\cdot h^{(m)} \\
&+ \phi\nabla(x_l-g_l^{(m)})\cdot h^{(m)}+\frac{\p\phi^{(m)}}{\p t}\, d\|V_t^{(m)}\|dt=:I^{(m)}_1+I_2^{(m)}+I_3^{(m)}+I_4^{(m)}+I_5^{(m)}.
\end{split}
\label{b22}
\end{equation}
We subsequently identify $\lim_{m\rightarrow\infty}(\mu^{(m)})^{-1}I_j^{(m)}$ for each $j=1,\cdots,5$. 
\newline
{\it Estimate of $I_1^{(m)}$}. 
\newline
By \eqref{b3} and \eqref{b13}, we have $|u^{(m)}|=o(\mu^{(m)})$. Moreover, by \eqref{b16} and \eqref{b20}, we have 
$|\phi^{(m)}|\leq 3c_4\mu^{(m)}\sup |\phi|$.  
Thus we have $|u^{(m)}|^2|\phi^{(m)}|=o((\mu^{(m)})^3)$ and since $\|V_t^{(m)}\|(C(T,1/2))$ is uniformly bounded, we have
\begin{equation}
\lim_{m\rightarrow\infty}(\mu^{(m)})^{-1}I_1^{(m)}=0.
\label{b23}
\end{equation}
\newline
{\it Estimate of $I_2^{(m)}$}.
\newline
By \eqref{b8} and \eqref{b19}, one observes that $|\nabla\phi^{(m)}|\leq O(1)$ on
${\rm spt}\, \|V^{(m)}_t\|$. Since $|u^{(m)}|=o(\mu^{(m)})$, we conclude that 
\begin{equation}
\lim_{m\rightarrow\infty}(\mu^{(m)})^{-1} I_2^{(m)}=0.
\label{b24}
\end{equation}
\newline
{\it Estimate of $I_3^{(m)}$}.
\newline
Choose $\tilde{\phi}\in C_c^{\infty}(B_{1/2}^k;\R^+)$ such that $\tilde{\phi}=1$ on ${\rm spt}\, \phi(\cdot,t)$
for all $t\in(-3/4,1)$. 
We also re-define $\tilde{\phi}(x,t):=\tilde{\phi}(T(x))$ for $x\in C(T,1/2)$. Take $-3/4<t_1<t_2<1$ so that
${\rm spt}\, \phi\subset B_{1/2}^k\times (t_1,t_2)$. 
We first claim that
\begin{equation}
\lim_{m\rightarrow\infty}\int_{t_1}^{t_2}\int_{C(T,1/2)}\tilde{\phi}|h^{(m)}|^2\, d\|V_t^{(m)}\|dt=0.
\label{b25}
\end{equation}
For the proof, using $\tilde{\phi}$ in \eqref{main}, we have
\begin{equation}
\left.\int_{C(T,1/2)}\tilde{\phi}\, d\|V_t^{(m)}\|\right|_{t=t_1}^{t_2}\leq \int_{t_1}^{t_2}\int_{C(T,1/2)}(-h^{(m)}\tilde{\phi}+\nabla\tilde{\phi})\cdot(h^{(m)}+(u^{(m)})^{\perp})\, d\|V_t^{(m)}\|dt.
\label{b26}
\end{equation}
By \eqref{b7}, we have $d\|V^{(m)}_t\|\rightarrow d{\mathcal H}^k
\lfloor_T$  uniformly in $t$ and the left-hand side of \eqref{b26}
converges to 0 as $m\rightarrow \infty$. The right-hand side of \eqref{b26} is bounded from above by
\begin{equation}
\begin{split}
&\int_{t_1}^{t_2}\int_{C(T,1/2)}-|h^{(m)}|^2\tilde{\phi}+|u^{(m)}|^2\tilde{\phi}+\frac14|h^{(m)}|^2\tilde{\phi}+|u^{(m)}||\nabla\tilde{\phi}|
+|h^{(m)}||(\nabla\tilde{\phi})^{\perp}|\, d\|V_t^{(m)}\|dt \\
& \leq \int_{t_1}^{t_2}\int_{C(T,1/2)}-\frac12 |h^{(m)}|^2\tilde{\phi}+|u^{(m)}|^2\tilde{\phi}+|u^{(m)}||\nabla\tilde{\phi}|+
\frac{|(\nabla\tilde{\phi})^{\perp}|^2}{\tilde{\phi}}\, d\|V_t^{(m)}\|dt,
\end{split}
\label{b27}
\end{equation}
where we used \eqref{perpH}. Terms involving $u^{(m)}$ converge to 0 since $|u^{(m)}|=o(\mu^{(m)})$. We also have
$(\nabla\tilde{\phi})^{\perp}=(T-{\rm image}\,\nabla f^{(m)})(\nabla\tilde{\phi})$ since $\nabla\tilde{\phi}=T(\nabla\tilde{\phi})$. By \eqref{b7}, we have
\begin{equation}
\frac{|(\nabla\tilde{\phi})^{\perp}|^2}{\tilde{\phi}}\leq \|T-{\rm image}\,\nabla f^{(m)}\|^2 \frac{|\nabla\tilde{\phi}|^2}{\tilde{\phi}}\rightarrow
0
\label{b28}
\end{equation}
as $m\rightarrow\infty$ uniformly in $(x,t)$. Combining \eqref{b26}-\eqref{b28}, we prove \eqref{b25}. 
Since we took $\tilde{\phi}$ so that ${\rm spt}\, \phi\subset \{\tilde{\phi}=1\}$, \eqref{b16} and \eqref{b25} show that
\begin{equation}
\lim_{m\rightarrow\infty} (\mu^{(m)})^{-1} I_3^{(m)}=0.
\label{b29}
\end{equation}
\newline
{\it Estimate of $I_{4}^{(m)}+I_5^{(m)}$}.
\newline
By arguing via the Gram-Schmidt orthonormalization as in the proof of Lemma \ref{propQ}, 
one can show that there exists a constant $c_5=c(n,k)$ such that
\begin{equation}
|{\rm image}\, \nabla f^{(m)}-\sum_{j=1}^k {\bf f}_j^{(m)}\otimes {\bf f}_j^{(m)} |\leq c_5 |\nabla f^{(m)}|^2,
\label{b30}
\end{equation}
where ${\bf f}_j^{(m)}=(0,\cdots,1,\cdots,0,\frac{\p f^{(m)}_{k+1}}{\p x_j},\cdots,\frac{\p f^{(m)}_n}{\p x_j})$ with $1$ in the
$j$-th component. Here we recall that we are identifying ${\rm image}\, \nabla f^{(m)}$ with the 
corresponding $n\times n$ orthogonal projection matrix. Then we use \eqref{first} to derive
\begin{equation}
I_4^{(m)}=\int_{-3/4}^1\int_{C(T,1/2)}(-\nabla\phi\otimes \nabla (x_l-g_l^{(m)}) +\phi\nabla^2 g_l^{(m)})\cdot 
({\rm image}\, \nabla f^{(m)})
\, d\|V_t^{(m)}\|dt.
\label{b31}
\end{equation}
We also have
\begin{equation}
-\nabla\phi\otimes\nabla(x_l-g_l^{(m)})\cdot\sum_{j=1}^k {\bf f}_j^{(m)}\otimes {\bf f}_j^{(m)}=
\nabla( g_l^{(m)}-f_l^{(m)})\cdot \nabla \phi
\label{b32}
\end{equation}
and
\begin{equation}
\phi\nabla^2 g_l^{(m)}\cdot \sum_{j=1}^k {\bf f}_j^{(m)}\otimes
{\bf f}_j^{(m)}=\phi \Delta g_l^{(m)}.
\label{b33}
\end{equation}
Since $g_l^{(m)}$ satisfies \eqref{propfl2}, we obtain by \eqref{b30}-\eqref{b33} and \eqref{b8} that
\begin{equation}
\begin{split}
&\left|I_4^{(m)}+I_5^{(m)}-\int_{-3/4}^1\int_{C(T,1/2)}\nabla(g_l^{(m)}-f_l^{(m)})\cdot
\nabla \phi +\mu^{(m)}(\tilde{f}_l^{(m)}+2c_4)\frac{\p \phi}{\p t}\, d\|V_t^{(m)}\|dt\right| \\
&\leq c_5 (\|\phi\|_{0}+\|\nabla\phi\|_{0}+\|\frac{\p \phi}{\p t}\|_{0})\|\nabla f^{(m)}\|^2_{0}.
\end{split}
\label{b34}
\end{equation}
By \eqref{b13}, \eqref{b17} and \eqref{b34}, we obtain
\begin{equation}
\lim_{m\rightarrow\infty}(\mu^{(m)})^{-1}(I_4^{(m)}+I_5^{(m)})=\int_{-3/4}^1\int_{B_{1/2}^k}
-\nabla \tilde{f}_l\cdot \nabla\phi+(\tilde{f}_l+2c_4)\frac{\p \phi}{\p t}\, d{\mathcal H}^k dt.
\label{b35}
\end{equation}
By \eqref{b22}, \eqref{b23}, \eqref{b24}, \eqref{b29} and \eqref{b35} and noting that
$\phi$ has a compact support in $B_{1/2}^k\times(-3/4,1)$, we obtain
\begin{equation}
0\leq \int_{-3/4}^1\int_{B_{1/2}^k}
-\nabla \tilde{f}_l\cdot \nabla\phi+\tilde{f}_l\frac{\p \phi}{\p t}\, d{\mathcal H}^k dt.
\label{b36}
\end{equation}
We may repeat the same computation with $\phi^{(m)}$ in \eqref{b19} replaced by
$(g_l^{(m)}-x_l+2c_4\mu^{(m)})\phi$, which is again nonnegative on ${\rm spt}\, \|V_t^{(m)}\|$. This leads
to the same conclusion as in \eqref{b36} with $\tilde{f}_l$ there replaced by $-\tilde{f}_l$.
Thus \eqref{b36} holds with equality for all $\phi\in C^{\infty}_c(B_{1/2}^k\times(-3/4,1);{\R}^+)$. 
This proves that $\tilde{f}_l$ satisfies the
heat equation in a weak sense. By the standard parabolic regularity theory, $\tilde{f}_l$
is $C^{\infty}$ and is a classical solution. This concludes the proof of Lemma \ref{heateq}.
\hfill{$\Box$}

We next prove 
\begin{lemma}
\begin{equation}
\lim_{m\rightarrow\infty} \|\tilde{f}^{(m)}-\tilde{f}\|_{L^2(B_{1/2}^k\times(-3/4,1))}=0.
\label{b37}
\end{equation}
\label{strong}
\end{lemma}
{\it Proof of Lemma \ref{strong}}. Let $l\in \{k+1,\cdots,n\}$ be fixed. 
We first claim that for each $s\in (-3/4,1)$, 
\begin{equation}
\tilde{f}^{(m)}_l(\cdot,s) \rightharpoonup \tilde{f}_l(\cdot,s)\mbox{ weakly in }
L^2(B_{1/2}^k).
\label{b38}
\end{equation}
By \eqref{b16}, for each fixed $s\in (-3/4,1)$, $\{\tilde{f}_l^{(m)}(\cdot,s)\}_{m=1}^{\infty}$ is bounded
in $L^2(B_{1/2}^k)$ in particular. Let $w\in L^2(B_{1/2}^k)$ be any weak limit. Let $\phi^{(m)}$
be as in \eqref{b20} and use \eqref{main} for $t_1=-3/4$ and $t_2=s$. The same computations
\eqref{b21}-\eqref{b35} show that we have
\begin{equation}
\int_{B_{1/2}^k} (w+2c_4)\phi(\cdot,s)\, d{\mathcal H}^k\leq \int_{-3/4}^s\int_{B_{1/2}^k}
-\nabla\tilde{f}_l\cdot\nabla\phi+(\tilde{f}_l+2c_4)\frac{\p\phi}{\p t}\, d{\mathcal H}^k dt.
\label{b39}
\end{equation}
Since $\tilde{f}_l$ is already known to be the solution of the heat equation, we have from \eqref{b39}
\begin{equation}
\int_{B_{1/2}^k}(w+2c_4)\phi(\cdot,s)\, d{\mathcal H}^k\leq \int_{B_{1/2}^k}(\tilde{f}_l(\cdot,s)+2c_4)\phi(\cdot,s)\,
d{\mathcal H}^k.
\label{b40}
\end{equation}
Similarly, replacing $\phi^{(m)}$ in \eqref{b19} by $(g_l^{(m)}-x_l+2c_4\mu^{(m)})\phi$, we obtain 
\begin{equation}
\int_{B_{1/2}^k}(2c_4-w)\phi(\cdot,s)\, d{\mathcal H}^k\leq \int_{B_{1/2}^k}(2c_4-\tilde{f}_l(\cdot,s))\phi(\cdot,s)\,
d{\mathcal H}^k.
\label{b41}
\end{equation}
Thus \eqref{b40} and \eqref{b41} show that $\int_{B_{1/2}^k}(w-\tilde{f}_l(\cdot,s))\phi(\cdot,s)\, d{\mathcal H}^k=0$. Since $\phi(\cdot,s)\in C^1_c(B_{1/2}^k;{\R}^+)$ may be 
chosen arbitrarily, we proved $w=\tilde{f}_l(\cdot,s)$ a.e. on $B_{1/2}^k$. Since any weak subsequence converges to $\tilde{f}_l(\cdot,s)$, the
whole sequence converges weakly to $\tilde{f}_l(\cdot,s)$, proving \eqref{b38}. The lower semicontinuity under
weak convergence shows that
\begin{equation}
\|\tilde{f}_l(\cdot,s)\|_{L^2(B_{1/2}^k)}\leq \liminf_{m\rightarrow\infty}\|\tilde{f}^{(m)}_l(\cdot,s)\|_{L^2(B_{1/2}^k)}
\label{b42}
\end{equation}
for all $s\in (-3/4,1)$. 
We next show that for any $-3/4<s_j<s<1$ with $s_j$ satisfying \eqref{b17.5} and for any $\phi\in
C^{\infty}_c (B_{1/2}^k;{\R}^+)$, we have
\begin{equation}
\limsup_{m\rightarrow\infty}\|\phi\tilde{f}_l^{(m)}(\cdot,s)\|^2_{L^2(B_{1/2}^k)}\leq 
\|\phi\tilde{f}_l(\cdot,s_j)\|^2_{L^2(B_{1/2}^k)}+c(\phi)(s-s_j).
\label{b43}
\end{equation}
To prove \eqref{b43}, we use $(x_l-g_l^{(m)})^2\phi$ as a test function in \eqref{main} with time interval $[s_j,s]$. Then by
the Cauchy-Schwarz inequality and dropping the positive term, we obtain
\begin{equation}
\begin{split}
&\left.\int_{C(T,1/2)}(x_l-g_l^{(m)})^2\phi\, d\|V_t^{(m)}\|\right|_{t=s_j}^s  \leq \int_{s_j}^s\int_{C(T,1/2)}
(x_l-g_l^{(m)})^2\phi|u^{(m)}|^2 \\ &+|u^{(m)}||\nabla((x_l-g_l^{(m)})^2\phi)| 
 +h^{(m)}\cdot \nabla ((x_l-g_l^{(m)})^2\phi)+\phi\frac{\p}{\p t}(x_l-g_l^{(m)})^2\, d\|V_t^{(m)}\|dt.
\end{split}
\label{b44}
\end{equation}
We divide both sides of \eqref{b44} by $(\mu^{(m)})^2$ and take $m\rightarrow\infty$. By \eqref{b38}, 
\eqref{b17.5} and \eqref{b7}, we have
\begin{equation}
\left.\int_{B_{1/2}^k}(\tilde{f}_l(\cdot,t))^2\phi\, d{\mathcal H}^k\right|_{t=s_j}^s\leq \liminf_{m\rightarrow\infty}(\mu^{(m)})^{-2}
\left.\int_{C(T,1/2)}(x_l-g_l^{(m)})^2\phi\, d\|V_t^{(m)}\|\right|_{t=s_j}^s
\label{b45}
\end{equation}
where we emphasize that the strong convergence at $t=s_j$ is essentially used. By \eqref{b16},
\eqref{b3}, \eqref{b13} and \eqref{b8}, one can check that 
\begin{equation}
\lim_{m\rightarrow\infty}(\mu^{(m)})^{-2}\int_{C(T,1/2)}(x_l-g_l^{(m)})^2\phi|u^{(m)}|^2
+|u^{(m)}||\nabla((x_l-g_l^{(m)})^2\phi)|\, d\|V_t^{(m)}\|=0
\label{b46}
\end{equation}
uniformly in $t$. 
For the last two terms of \eqref{b44}, by \eqref{first}, we have for a.e. $t\in (s_j,s)$
\begin{equation}
\begin{split}
& \int_{C(T,1/2)} h^{(m)}\cdot\nabla((x_l-g_l^{(m)})^2\phi)+\phi\frac{\p}{\p t}(x_l-g_l^{(m)})^2\,
d\|V_t^{(m)}\| \\
&=\int_{C(T,1/2)}-({\rm image}\, \nabla f^{(m)})\cdot \nabla^2 ((x_l-g_l^{(m)})^2\phi)
+\phi\frac{\p}{\p t}(x_l-g_l^{(m)})^2\,d\|V_t^{(m)}\|.
\end{split}
\label{b47}
\end{equation}
For any $S\in {\bf G}(n,k)$, 
\begin{equation}
\begin{split}
-S\cdot & \nabla^2 ((x_l-g_l^{(m)})^2\phi)
\leq -2|S(\nabla(x_l-g_l^{(m)}))|^2\phi+2(x_l-g_l^{(m)})\phi S\cdot\nabla^2 g_l^{(m)} \\
 &-(x_l-g_l^{(m)})^2 S\cdot\nabla^2\phi+2|S(\nabla(x_l-g_l^{(m)}))|^2\phi+2(x_l-g_l^{(m)})^2
\frac{|\nabla\phi|^2}{\phi}\\
& =2(x_l-g_l^{(m)})\phi S\cdot\nabla^2 g_l^{(m)}
+(x_l-g_l^{(m)})^2(-S\cdot\nabla^2\phi+2\frac{|\nabla\phi|^2}{\phi}).
\end{split}
\label{b48}
\end{equation}
For $S={\rm image}\, \nabla f^{(m)}$, by Lemma \ref{propQ}, we have 
\begin{equation}
|(S-T)\cdot\nabla^2 g^{(m)}_l|\leq c(n,k)
|\nabla f^{(m)}|^2 |\nabla^2 g^{(m)}_l|.
\label{b49}
\end{equation}
Thus \eqref{b48} and \eqref{b49} show that
\begin{equation}
\begin{split}
-({\rm image}\, \nabla f^{(m)}) & \cdot\nabla^2 ((x_l-g_l^{(m)})^2\phi)
\leq 2(x_l-g_l^{(m)})\phi T\cdot\nabla^2 g_l^{(m)} \\ &+c_6(n,k,\|\phi\|_{C^2})
\{|x_l-g_l^{(m)}||\nabla f^{(m)}|^2|\nabla^2 g_l^{(m)}|+(x_l-g_l^{(m)})^2\}.
\end{split}
\label{b50}
\end{equation}
By \eqref{propfl2}, we have $T\cdot\nabla^2 g_l^{(m)}-\frac{\p g_l^{(m)}}{\p t}=0$, thus
\eqref{b47}, \eqref{b50}, \eqref{b13}, \eqref{b8} and \eqref{b16} show that
\begin{equation}
\begin{split}
\limsup_{m\rightarrow\infty}   (\mu^{(m)})^{-2}\int_{s_j}^s\int_{C(T,1/2)} &  h^{(m)}\cdot\nabla((x_l-g_l^{(m)})^2\phi)+\phi\frac{\p}{\p t}(x_l-g_l^{(m)})^2\,
d\|V_t^{(m)}\| dt \\ 
\leq c_6 (c_4)^2 (s-s_j).
\end{split}
\label{b51}
\end{equation}
By combining \eqref{b44}, \eqref{b45}, \eqref{b46} and \eqref{b51}, we obtain \eqref{b43}. 
Using the smoothness of $\tilde{f}_l$, \eqref{b42}, \eqref{b43} and the fact that $\{s_j\}_{j=1}^{\infty}$
is dense, one can prove 
$\|\tilde{f}^{(m)}_l(\cdot,s)\|_{L^2(B_{1/2}^k)}\rightarrow\|\tilde{f}_l(\cdot,s)\|_{L^2(B_{1/2}^k)}$
for all $s\in (-3/4,1)$, which shows the strong $L^2(B_{1/2}^k)$ convergence. Since these $L^2$ norms are all
bounded uniformly in $s$ by \eqref{b16}, the dominated convergence theorem proves
the desired strong convergence, \eqref{b37}. This concludes the proof of Lemma \ref{strong}.
\hfill{$\Box$}

Next define $p=(p_{k+1},\cdots,p_n):{\R}^k\times\R\rightarrow \R^{n-k}$ by 
\begin{equation}
p_l(x,t):=\tilde{f}_l(0,0)+\nabla\tilde{f}_l(0,0)\cdot x+\frac12 \nabla^2\tilde{f}_l(0,0)\cdot x\otimes x+\frac{\p \tilde{f}_l}{\p t}(0,0)t
\label{b52}
\end{equation}
for $l=k+1,\cdots,n$. Since $\tilde{f}_l$ satisfies the heat equation by Lemma \ref{heateq}, the standard interior 
estimates with \eqref{b17.1} gives
\begin{equation}
\|p\|_{C^2(B_1^k\times(-1,1))}\leq c_7(n,k).
\label{b53}
\end{equation}
By the Taylor theorem again with the standard interior estimates, we have for $0<\theta<1/4$
\begin{equation}
\sup_{B_{\theta}^k\times(-\theta^2,\theta^2)}|\tilde{f}-p|\leq c_8(n,k)\theta^3.
\label{b54}
\end{equation}
We define for each $m\in \N$
\begin{equation}
\hat{g}^{(m)}:=g^{(m)}+\mu^{(m)}p
\label{b55}
\end{equation}
and define $Q_{\hat{g}^{(m)}}$ as in \eqref{defQ}. By \eqref{b52} and Lemma \ref{heateq},
we have $\hat{g}^{(m)}\in \F$. On ${\rm spt}\, \|V_t^{(m)}\|$, by \eqref{b55} and \eqref{b15}, 
\begin{equation}
Q_{\hat{g}^{(m)}}=\frac12 |f^{(m)}-g^{(m)}-\mu^{(m)}p|^2\leq (\mu^{(m)})^2(|\tilde{f}^{(m)}-\tilde{f}|^2+|\tilde{f}-p|^2).
\label{b56}
\end{equation}
Thus Lemma \ref{strong}, \eqref{b7}, \eqref{b54} and \eqref{b56} show
\begin{equation}
\limsup_{m\rightarrow\infty}(\mu^{(m)})^{-2}\int_{-\theta^2}^{\theta^2}\int_{C(T,\theta)} Q_{\hat{g}^{(m)}}\, d\|V_t^{(m)}\|dt
\leq 2\omega_k (c_8)^2 \theta^{k+8}.
\label{b57}
\end{equation}
We now choose a small $0<\theta<1/4$ depending only on $n,\, k,\, \alpha$ so that
\begin{equation}
1\geq 4\omega_k(c_8)^2 \theta^{2-2\alpha}
\label{b58}
\end{equation}
holds. By \eqref{b53} and \eqref{b55}, we have $\|\hat{g}^{(m)}-g^{(m)}\|_{C^2(B_1^k\times(-1,1))}
\leq c_7\mu^{(m)}$, thus \eqref{b10} is satisfied for all sufficiently large $m$. One can check 
that \eqref{b11} and \eqref{b57} lead to a contradiction due to \eqref{b58} for all sufficiently
large $m$. Thus we complete the proof of Proposition \ref{blow}.
\hfill{$\Box$}
\section{$C^{2,\alpha}$ estimate}
Working under the same conditions as in the previous section and iterating the 
argument, we show a proper decay properties necessary for the proof of $C^{2,\a}$
estimates. First we prove
\begin{prop}
Corresponding to $n,\, k,\,\alpha$ there exist $0<\varepsilon_2<1$ and $1<c_{9}<\infty$ 
with the following property. Under the assumptions of Proposition \ref{blow} with $\varepsilon_1$
replaced by $\varepsilon_2$, with additional assumptions
\begin{equation}
R^{1+\alpha}[u]_{\alpha}\leq \varepsilon_2,
\label{cb0}
\end{equation}
\begin{equation}
\nabla f(0,0)=0, 
\label{cb1}
\end{equation}
\begin{equation}
R^{1+\alpha}[\nabla f]^2_{\frac{(1+\alpha)}{2}}\leq \varepsilon_2,
\label{cb1.1}
\end{equation}
we have
\begin{itemize}
\item[(1)] at $(x,t)=(0,0)$, $f$ is differentiable with respect to $(x,t)$ and 
$\nabla f$ is differentiable with respect to $x$,
\item[(2)] there exists $g^{(0)}\in \F$ such that 
\begin{equation}
f=g^{(0)},\hspace{.3cm} \nabla g^{(0)}=0,\hspace{.3cm} \nabla^2 f=\nabla^2 g^{(0)},
\hspace{.3cm} \frac{\p f}{\p t}=\frac{\p g^{(0)}}{\p t}
\label{cb2}
\end{equation}
all hold at $(x,t)=(0,0)$ and 
\begin{equation}
\begin{split}
R^{-1}\|g-g^{(0)}\|_0+\|\nabla( g-g^{(0)})\|_0 & +R\|\nabla^2( g-g^{(0)})\|_0+R\|\frac{\p}{\p t}
(g-g^{(0)})\|_0 \\
&\leq c_9 \max\{\mu, c_3 R^{1+\alpha}[u]_{\alpha}, c_3 R^{1+\alpha}[\nabla f]_{\frac{(1+\alpha)}{2}}^2\}.
\end{split}
\label{cb3}
\end{equation}
Here $\|\cdot\|_0:=\sup_{B_R^k\times(-R^2,R^2)}|\cdot |$. 
\item[(3)] Whenever $0<r\leq R$, there exists $g^{(r)}\in \F$ such that 
\begin{equation}
\begin{split}
 r^{-1}\|g^{(r)}-g^{(0)}\|_0&+\|\nabla (g^{(r)}- g^{(0)})\|_0+r\|\nabla^2 (g^{(r)}- g^{(0)})\|_0
+r \|\frac{\p}{\p t}(g^{(r)}-g^{(0)})\|_0\\
&+\Big(r^{-k-4}\int_{-r^2}^{r^2}\int_{C(T,r)} Q_{g^{(r)}}\, d\|V_t\|dt\Big)^{\frac12} \\
&\leq c_9 (r/R)^{1+\alpha}\max\{\mu,c_3 R^{1+\alpha}[u]_{\alpha},c_3 R^{1+\alpha}[\nabla f]_{\frac{(1+\alpha)}{2}}^2\}.
\end{split}
\label{cb4}
\end{equation}
Here $\|\cdot\|_0:=\sup_{B_r^k\times(-r^2,r^2)}|\cdot|$. 
\end{itemize}
\label{longprop}
\end{prop}
{\it Proof}. After a change of variables, we may assume that $R=1$. For any $q\in C^2$ define
\begin{equation}
\| q\|_{C^{2,1}(r)}:=r^{-1}\| q\|_0+\|\nabla q\|_0
+r\|\nabla^2 q\|_0+r\|\frac{\p q}{\p t}\|_0
\label{cb5}
\end{equation}
with $\|\cdot\|_0=\sup_{B_r^k\times(-r^2,r^2)}|\cdot|$ here. For notational simplicity define
\begin{equation}
K:=\max\{\mu,c_3[u]_{\alpha},
c_3[\nabla f]_{\frac{(1+\alpha)}{2}}^2\}.
\label{K}
\end{equation}
We choose $0<\varepsilon_2<\varepsilon_1$ and $1<c_9<\infty$ so that
\begin{equation}
c_{3}\varepsilon_2\leq \varepsilon_1,
\label{cb5.1}
\end{equation}
\begin{equation}
\varepsilon_2+(c_3)^2 \varepsilon_2 \sum_{j=1}^{\infty}\theta^{(j-1)(1+\alpha)}\leq \varepsilon_1,
\label{cb5.2}
\end{equation}
\begin{equation}
6c_3\theta^{-2}\sum_{j=0}^{\infty}\theta^{(j-1)\alpha}\leq c_9,
\label{cb5.3}
\end{equation}
\begin{equation}
2\theta^{-(\frac{k}{2}+3+\alpha)}\leq c_9.
\label{cb5.4}
\end{equation}
We inductively prove the following claims. We set $g^{(1)}:=g$ and suppose that for $j=
1,\cdots,m$, there are $g^{(\theta^j)}\in \F$ such that 
\begin{equation}
\|g^{(\theta^{j})}-g^{(\theta^{j-1})}\|_{C^{2,1}(\theta^{j-1})}\leq c_3 \theta^{(j-1)(1+\alpha)}K,
\label{cb6}
\end{equation}
\begin{equation}
\mu^{(j)}:=\Big(\theta^{-j(k+4)}\int_{-\theta^{2j}}^{\theta^{2j}}\int_{C(T,\theta^j)} Q_{g^{(\theta^j)}}\, d\| V_t\|dt\Big)^{\frac12}
\leq \theta^{j(1+\alpha)}K. 
\label{cb7}
\end{equation}
Consider the case $j=1$. 
Since $\varepsilon_2\leq \varepsilon_1$, Proposition \ref{blow} gives $\hat{g}\in \F$ which we denote
by $g^{(\theta)}$. Note here that
due to \eqref{cb1}, we have $\|\nabla f\|_{0}\leq [\nabla f]_{\frac{(1+\alpha)}{2}}$. 
Assume \eqref{cb6} and \eqref{cb7} hold up to  $j=m$.  We have
\begin{equation}
\mu^{(m)}\leq \theta^{m(1+\alpha)}K\leq \varepsilon_1
\label{cb8}
\end{equation}
by \eqref{b4}, \eqref{cb0}, \eqref{cb1.1} and \eqref{cb5.1}. With the 
notation $\|\cdot\|=\sup_{B_{\theta^m}^k\times(-\theta^m,\theta^m)}|\cdot|$ for abbreviation
in the next computations, we compute
\begin{equation}
\begin{split}
&\|\nabla g^{(\theta^m)}\|+\theta^m\|\nabla^2 g^{(\theta^m)}\|+\theta^m\|\frac{\p g^{(\theta^m)}}{\p t}\| \\
&\leq \sum_{j=1}^{m}\big(\|\nabla(g^{(\theta^j)}-g^{(\theta^{j-1})})\|+\theta^m \|\nabla^2(g^{(\theta^j)}
-g^{(\theta^{j-1})})\|+\theta^m\|\frac{\p}{\p t}(g^{(\theta^j)}-g^{(\theta^{j-1})})\|\big)\\
&+\|\nabla g^{(1)}\|+\theta^m \|\nabla^2 g^{(1)}\|+\theta^m \|\frac{\p g^{(1)}}{\p t} \| \\
& \leq \varepsilon_2+
\sum_{j=1}^m c_3\theta^{(j-1)(1+\alpha)}K
\mbox{ (by \eqref{cb6} and \eqref{b2})}\\
& \leq \varepsilon_2+ (c_3)^2 \varepsilon_2\sum_{j=1}^{\infty}\theta^{(j-1)(1+\alpha)}
\leq \varepsilon_1\mbox{ (by \eqref{b4}, \eqref{cb0}, \eqref{cb1.1} and \eqref{cb5.2})}. 
\end{split}
\label{cb9}
\end{equation}
Hence for $R=\theta^m$, \eqref{b1}-\eqref{b4} are all satisfied due to \eqref{cb8} and \eqref{cb9}. 
By Proposition \ref{blow} there exists a new $\hat{g}\in \F$ denoted by $g^{(\theta^{m+1})}$ with the estimates
\begin{equation}
\|g^{(\theta^{m+1})}-g^{(\theta^m)}\|_{C^{2,1}(\theta^m)}\leq c_3 \mu^{(m)}\leq c_3
\theta^{m(1+\alpha)}K
\label{cb10}
\end{equation}
by \eqref{b5}, \eqref{cb8} and 
\begin{equation}
\begin{split}
\mu^{(m+1)} & \leq \theta^{1+\alpha}\max \{\mu^{(m)},c_3 \theta^{m(1+\alpha)}[u]_{\alpha},
c_3\sup_{B_{\theta^m}^k\times(-\theta^m,\theta^m)}|\nabla f|^2\} \\
& \leq \theta^{(m+1)(1+\alpha)}K
\end{split}
\label{cb11}
\end{equation}
by \eqref{b6}, \eqref{cb1} and \eqref{cb8}. \eqref{cb10} and \eqref{cb11} show that 
\eqref{cb6} and \eqref{cb7} are satisfied for $j=m+1$, thus they are satisfied for all 
$j\in \N$. For function $g$ satisfying \eqref{propfl1}, we have for $0< r\leq 1$ 
\begin{equation}
\|g\|_{C^{2,1}(1)}\leq |g(0,0)|+2|\nabla g(0,0)|+3 |\nabla^2 g(0,0)|+2|\frac{\p g}{\p t}(0,0)|
\leq 3 r^{-1} \|g\|_{C^{2,1}(r)}.
\label{cb12}
\end{equation}
Thus for $j\in \N$ we have by \eqref{cb6} and \eqref{cb12} 
\begin{equation}
\|g^{(\theta^{j})}-g^{(\theta^{j-1})}\|_{C^{2,1}(1)}\leq 3c_3 \theta^{(j-1)\alpha}K.
\label{cb13}
\end{equation}
It is clear from \eqref{cb13} that there exists $g^{(0)}=\lim_{j\rightarrow\infty}g^{(\theta^j)}$ which 
belongs to $\F$, and which satisfies \eqref{cb3} by \eqref{cb5.3}. For $0<r\leq 1$, choose $m$ such that
\begin{equation}
\theta^{m+1}< r\leq \theta^m
\label{cb13.5}
\end{equation}
and set $g^{(r)}=g^{(\theta^m)}$. Then by the similar
computations as in \eqref{cb12}, we have 
\begin{equation}
\begin{split}
&\|g^{(0)}-g^{(r)}\|_{C^{2,1}(r)}  \leq \sum_{j=m+1}^{\infty}\|g^{(\theta^j)}-g^{(\theta^{j-1})}\|_{C^{2,1}(r)}  
\\ &\leq \sum_{j=m+1}^{\infty}3 \theta^{m-j}\|g^{(\theta^{j})}-g^{(\theta^{j-1})}\|_{C^{2,1}(\theta^{j-1})} 
\leq \sum_{j=0}^{\infty}3c_3 \theta^{(m+1)(1+\alpha)-2}\theta^{(j-1)\alpha}K \\ & \leq \frac{c_9}{2} r^{1+\alpha}K
\end{split}
\label{cb14}
\end{equation}
by \eqref{cb6}, \eqref{cb5.3} and \eqref{cb13.5}. By \eqref{cb13.5}, \eqref{cb7} and \eqref{cb5.4}, we also have
\begin{equation}
\begin{split}
&\big(r^{-k-4}\int_{-r^2}^{r^2}\int_{C(T,r)}Q_{g^{(r)}}\, d\|V_t\|dt\big)^{\frac12}
  \leq \theta^{-(k+4)/2}\mu^{(m)} \\ &\leq \theta^{-(k/2+3+\alpha)}\theta^{(m+1)(1+\alpha)} K
\leq \frac{c_9}{2} r^{1+\alpha}K.
\end{split}
\label{cb15}
\end{equation}
Summation of \eqref{cb14} and \eqref{cb15} proves \eqref{cb4}. It is easy to check that \eqref{cb4} shows
$g^{(0)}(0,0)=f(0,0)$ and $\nabla g^{(0)}(0,0)=\nabla f(0,0)=0$. To prove the differentiability of $\nabla f$,
we will prove
\begin{equation}
\lim_{x\rightarrow 0}\frac{|\nabla f(x,0)-\nabla^2 g^{(0)}(0,0)\cdot x|}{|x|}=0,
\label{cb16}
\end{equation}
which also gives $\nabla^2 f(0,0)=\nabla^2 g^{(0)}(0,0)$. We note that 
$\nabla^2 g^{(0)}(0,0)\cdot x=\nabla g^{(0)}(x,0)$. By \eqref{cb4}, for $|x|<1/2$, 
$|\nabla g^{(0)}(x,0)-\nabla g^{(2|x|)}(x,0)|=O(|x|^{1+\alpha})$. Thus to show
\eqref{cb16}, it suffices to prove
\begin{equation}
\lim_{x\rightarrow0}\frac{|\nabla f(x,0)-\nabla g^{(2|x|)}(x,0)|}{|x|}=0.
\label{cb17}
\end{equation}
For any $x$ with $|x|<1/2$, set 
$r:=|x|$, $\beta:=\frac{2\alpha}{k+6}$ and let $A$ be the Affine $k$-dimensional plane 
which is tangent to the ${\rm graph}\, g^{(2r)}$ at $(x, g^{(2r)}(x,0))$. As a graph, $A$ is
represented as
\begin{equation}
z\in \R^k\rightarrow g^{(2r)}(x,0)+\nabla g^{(2r)}(x,0)\cdot (z-x).
\label{cb17.5}
\end{equation}
In the following we estimate
\begin{equation}
\mu(x):=\Big(r^{-(1+\beta)(k+4)}\int_{-r^{2(1+\beta)}}^{r^{2(1+\beta)}}
\int_{C(T,x,r^{1+\beta})}{\rm dist}\, (y,A)^2 \, d\|V_t\|(y)dt\Big)^{\frac12}
\label{cb17.6}
\end{equation}
to apply the gradient estimate of \cite[Th. 8.7]{Kasai-Tonegawa}. 
For $y=(z,f(z,t))\in {\rm spt}\, \|V_t\|$ with $z:=T(y)$, 
\begin{equation}
\begin{split}
{\rm dist}\, (y,A)& \leq  |f(z,t)-g^{(2r)}(x,0)-\nabla g^{(2r)}(x,0)\cdot (z-x)| \\
&\leq  |f(z,t)-g^{(2r)}(z,t)|+|g^{(2r)}(z,t)-g^{(2r)}(x,0)-\nabla g^{(2r)}(x,0)\cdot (z-x)| \\
& \leq \sqrt{2}Q_{g^{(2r)}}^{\frac12}(y,t)+\|\frac{\p g^{(2r)}}{\p t}\|_{0}|t|+\|\nabla^2 g^{(2r)}\|_{0}|z-x|^2 \\
&\leq \sqrt{2}Q_{g^{(2r)}}^{\frac12}(y,t) + c_{10}(|t|+|z-x|^2).
\end{split}
\label{cb18}
\end{equation}
The existence of $c_{10}$ independent of $r$ follows from \eqref{cb13}. Substituting \eqref{cb18} into 
\eqref{cb17.6} gives
\begin{equation}
\mu(x)\leq \Big(r^{-(1+\beta)(k+4)}\int_{-r^{2(1+\beta)}}^{r^{2(1+\beta)}}
\int_{C(T,x,r^{1+\beta})}2 Q_{g^{(2r)}}\, d\|V_t\|dt\Big)^{\frac12}+
6 c_{10}r^{1+\beta}.
\label{cb19}
\end{equation}
Since $C(T,x,r^{1+\beta})\subset C(T,2r)$, \eqref{cb19} shows
\begin{equation}
\mu(x)\leq 2 ^{\frac{k+5}{2}} r^{-\frac{\beta(k+4)}{2}}\Big((2r)^{-(k+4)}\int_{-4r^2}^{4r^2}
\int_{C(T,2r)} Q_{g^{(2r)}}\, d\| V_t\|dt\Big)^{\frac12}+6c_{10}r^{1+\beta}.
\label{cb20}
\end{equation}
By \eqref{cb4} and \eqref{cb20}, we obtain
\begin{equation}
\mu(x)\leq 2^{\frac{k+5}{2}}c_9 K r^{1+\alpha-\frac{\beta(k+4)}{2}}+6c_{10} r^{1+\beta}=(
2^{\frac{k+5}{2}}c_9 K+6c_{10}) r^{1+\beta}
\label{cb21}
\end{equation}
by the choice of $\beta=\frac{2\alpha}{k+6}$. Note that $u$ is H\"{o}lder continuous and
$u(0,0)=0$, thus for any large $p,\, q$ there which we subsequently fix, 
\begin{equation}
c_{11}(x):=(r^{1+\beta})^{1-\frac{k}{p}-\frac{2}{q}}\|u\|_{L^{p,q}(C(T,x, r^{1+\beta})\times (-r^{2(1+\beta)},r^{2(1+\beta)}))}
\leq \varepsilon_2 c(p,q,k,n) r^{1+\alpha}.
\label{cb22}
\end{equation}
The existence of $t_1$ and $t_2$ there for $\nu=3/4$, for example, is satisfied since 
$f$ is a graph with uniformly small spacial gradient. Thus there exists a constant $c_{12}$ depending only
on $n,\, k$ such that
\begin{equation}
|\nabla f(x,0)-\nabla g^{(2r)}(x,0)|\leq c_{12}\max\{\mu(x),c_{11}(x)\}
\leq 2c_{12} (
2^{\frac{k+5}{2}}c_9 K+6c_{10})r^{1+\beta}
\label{cb23}
\end{equation}
by \cite[Th. 8.7]{Kasai-Tonegawa}, \eqref{cb21} and \eqref{cb22} for all sufficiently small $r$.  
Now \eqref{cb23} proves \eqref{cb17}. Finally, we need to prove
\begin{equation}
\lim_{(x,t)\rightarrow 0}\frac{|f(x,t)-g^{(0)}(0,t)|}{\sqrt{|x|^2+t^2}}=0
\label{cb24}
\end{equation}
which will prove $f$ is differentiable at $(x,t)=(0,0)$ (recall $\nabla g^{(0)}(0,0)=0$) 
and $\frac{\p f}{\p t}(0,0)=\frac{\p g^{(0)}}{\p t}(0,0)$.
 Set $r:=(|x|^2+t^2)^{1/4}$. By \eqref{cb4}, we have for some $\hat{t}$ with $|\hat{t}|\leq r^2$ that
\begin{equation}
\begin{split}
|g^{(0)}(0,t)-g^{(2r)}(0,t)| & \leq |g^{(0)}(0,0)-g^{(2r)}(0,0)|+|t| \Big|\frac{\p}{\p t} (g^{(0)}-g^{(2r)})(0,\hat{t})\Big| \\
& \leq c_9 (2r)^{2+\alpha}K+2^{\alpha}c_9 r^{2+\alpha} K.
\end{split}
\label{cb25}
\end{equation}
Moreover, for some $\hat{x}$ with $|\hat{x}|\leq r^2$, 
\begin{equation}
|g^{(2r)}(x,t)-g^{(2r)}(0,t)|\leq |\nabla g^{(2r)}(\hat{x},t)||x|\leq c_9 2^{1+\alpha} r^{3+\alpha}K
\label{cb25.5}
\end{equation}
by \eqref{cb4}. 
Thus \eqref{cb25} and \eqref{cb25.5} show that it suffices to prove
\begin{equation}
\lim_{(x,t)\rightarrow0}\frac{|f(x,t)-g^{(2r)}(x,t)|}{\sqrt{|x|^2+t^2}}=0
\label{cb26}
\end{equation}
to prove \eqref{cb24}. 
We basically repeat the same argument as for $\nabla^2 f$. Set 
$\beta$ as before and let $A$ be the Affine $k$-dimensional plane
which is tangent to the ${\rm graph}\, g^{(2r)}$ at $(x, g^{(2r)}(x,t))$. 
We define
\begin{equation}
\mu(x,t):=\Big(r^{-(1+\beta)(k+4)}\int_{t-r^{2(1+\beta)}}^{t+r^{2(1+\beta)}}
\int_{C(T,x,r^{1+\beta})}{\rm dist}\, (y,A)^2\, d\|V_s\|ds\Big)^{\frac12}.
\label{cb27}
\end{equation}
For $y=(z,g(z,s))\in {\rm spt}\, \|V_s\|$ with $z:=T(y)$, 
\begin{equation}
\begin{split}
{\rm dist}\, (y,A) & \leq |f(z,s)-g^{(2r)}(x,t)-\nabla g^{(2r)}(x,t)\cdot
(z-x)| \\
&\leq |f(z,s)-g^{(2r)}(z,s)|+|g^{(2r)}(z,s)-g^{(2r)}(x,t)-\nabla g^{(2r)}(x,t)\cdot (z-x)| \\
& \leq \sqrt{2}Q_{g^{(2r)}}^{\frac12}(z,s)+|s-t|\|\frac{\p g^{(2r)}}{\p s}\|_{0}
+|z-x|^2 \|\nabla ^2 g^{(2r)}\|_{0} \\
& \leq \sqrt{2}Q_{g^{(2r)}}^{\frac12}(z,s)+c_{10}(|s-t|+|z-x|^2).
\end{split}
\label{cb28}
\end{equation}
Substitute \eqref{cb28} into \eqref{cb27}, and proceed just as before. Note
that 
\begin{equation}
C(T,x,r^{1+\beta})\times(t-r^{2(1+\beta)},t+r^{2(1+\beta)})\subset C(T,2r)\times (-4r^2,
4r^2)
\label{cb28.5}
\end{equation}
due to $r=(|x|^2+t^2)^{1/4}$. 
Then we obtain by \eqref{cb27}-\eqref{cb28.5} (with an obvious modification for $c_{11}(x,t)$ and $\sup$ estimate
instead of gradient estimate of \cite[Th. 8.7]{Kasai-Tonegawa}) 
\begin{equation}
r^{-(1+\beta)}|f(x,t)-g^{(2r)}(x,t)|\leq c_{12}\max\{\mu(x,t),c_{11}(x,t)\}
\leq 2 c_{12}(2^{\frac{k+5}{2}} c_{9}K+6 c_{10}) r^{1+\beta}
\label{cb29}
\end{equation}
for all sufficiently small $r$. By \eqref{cb29}, we prove \eqref{cb26}.
This completes the proof of Proposition \ref{longprop}.
\hfill{$\Box$}

Next, to apply the estimates of Proposition \ref{longprop} at a given point, we need to make a 
change of variables so that $\nabla f$ and $u$ are both zero there with respect
to the new coordinate system. Suppose that we 
have $\{V_t\}_{-1<t<1}$ and $\{u(\cdot,t)\}_{-1<t<1}$ satisfying (B1)-(B4) on
$B_1\times(-1,1)$. Let $(\tilde{x},\tilde{t})\in B_{1/2}\cap{\rm spt}\, \|V_t\|$ with
$-1/2<\tilde{t}<1/2$ be arbitrary. 
By suitable rotation and parallel translation, we may 
choose a coordinate system so that $(\tilde{x},\tilde{t})$ is translated to
the origin $(0,0)$ and ${\rm spt}\, \|V_{0}\|$ is
tangent to $\R^k\times\{0\}$, so that the graph of $f$ has $\nabla f(0,0)=0$.
Note that $B_{1/2}\times(-1/2,1/2)$ (with this new coordinate system) is included in the
original domain. To have $u(0,0)=0$, we change 
the variables by $(x,t)\rightarrow (x-tu(0,0),t)$. Namely, we introduce a
new coordinate system so that the frame moves at the constant speed $u(0,0)$.
Define for each $t\in (-1/2,1/2)$ and $\phi\in C_c(G_k(\R^n))$
\begin{equation}
\tilde{V}_t (\phi):=V_t(\phi(\cdot-tu(0,0),\cdot)),
\hspace{.3cm}
\tilde{u}(x,t):=u(x+tu(0,0),t)-u(0,0).
\label{t1}
\end{equation}
If $u(0,0)$ is assumed to be sufficiently
small, $B_{1/4}\times(-1/4,1/4)$ is included in the original domain under the
new coordinate system. It is 
natural to expect the following.
\begin{lemma}
The newly defined $\{\tilde{V}_t\}_{-1/4<t<1/4}$ and $\{\tilde{u}(\cdot,t)\}_{-1/4<t<1/4}$ satisfy (B4)
on $B_{1/4}\times(-1/4,1/4)$ and $\tilde{u}(0,0)=0$. 
\label{coor}
\end{lemma}
{\it Proof}. Obviously $\tilde{u}(0,0)=0$ follows from \eqref{t1}. Write 
$a:=u(0,0)$ for simplicity. We need to check 
that \eqref{main} holds for $\tilde{V}_t$ and $\tilde{u}$. For any $\phi\in C^1(B_{1/4}
\times(-1/4,1/4);\R^+)$ with $\phi(\cdot, t)\in C^1_c(B_{1/4})$, 
define $\tilde{\phi}(x,t):=\phi(x-at,t)$. Then for any $-1/4<t_1<t_2<1/4$,
by \eqref{t1} and \eqref{main},
\begin{equation}
\left.\|\tilde{V}_{t}\|(\phi(\cdot,t))\right|_{t=t_1}^{t_2}=\left.\|V_t\|(\tilde{\phi}(\cdot,t))
\right|_{t=t_1}^{t_2}\leq \int_{t_1}^{t_2}\int
(\nabla\tilde{\phi}-\tilde{\phi}h)\cdot (h+u^{\perp})+\frac{\p \tilde{\phi}}{\p t}\, d\|V_t\|dt.
\label{t1.5}
\end{equation}
If we denote the mean curvature vector of $\tilde{V}_t$ by $\tilde{h}(\tilde{V}_t, x)$, 
we have $\tilde{h}(\tilde{V}_t,x-at)=h(V_t, x)$ since the change of variables is simply
a translation for each fixed time. Thus
\begin{equation}
\begin{split}
& \int(\nabla\tilde{\phi}-\tilde{\phi}h)\cdot (h+u^{\perp})\, d\|V_t\|\\
& =\int
\{\nabla \tilde{\phi}-\tilde{\phi}\tilde{h}(\tilde{V}_t,\cdot-at)\}
\cdot (\tilde{h}(\tilde{V}_t,\cdot-at)+\tilde{u}^{\perp}(\cdot-at)+a^{\perp})\, d\|V_t\| \\
&=\int (\nabla \phi-\phi\tilde{h})\cdot (\tilde{h}+\tilde{u}^{\perp}+a^{\perp})\, d\|\tilde{V}_t\|.
\end{split}
\label{t2}
\end{equation}
By \eqref{first} and \eqref{perpH} on the other hand, for a.e$.$ $t$, we have
\begin{equation}
\begin{split}
& \int(\nabla\phi-\phi\tilde{h})\cdot a^{\perp}\, d\|\tilde{V}_t\|=\int\nabla\phi\cdot a^{\perp}
-\phi\tilde{h}\cdot a\, d\|\tilde{V}_t\| \\ 
& =\int\nabla\phi\cdot a^{\perp}\, d\|\tilde{V}_t\|
+\int S\cdot (a\otimes\nabla \phi)\, d\tilde{V}_t(\cdot,S)
=\int\nabla\phi\cdot a\, d\|\tilde{V}_t\|
\end{split}
\label{t3}
\end{equation}
since $S(a\otimes \nabla \phi)=\nabla\phi \cdot a^{\top}$, where $\cdot^{\top}$ is the
projection to the tangent space, and $a^{\perp}+a^{\top}=a$. Since $\frac{\p\tilde{\phi}}{\p t}(x,t)=-\nabla\phi(\cdot-at,t)\cdot
a+\frac{\p\phi}{\p t}(\cdot-at,t)$, \eqref{t1.5}-\eqref{t3} prove
\begin{equation}
\left.\|\tilde{V}_t\|(\phi(\cdot,t))\right|_{t=t_1}^{t_2}\leq \int_{t_1}^{t_2}\int
(\nabla\phi-\phi\tilde{h})\cdot(\tilde{h}+\tilde{u}^{\perp})+\frac{\p\phi}{\p t}\, d\|\tilde{V}_t\|dt.
\label{t4}
\end{equation}
\eqref{t4} shows the claim of the present lemma. 
\hfill{$\Box$}

Finally, assuming that we already have $C^{1,\frac{1+\a}{2}}$ estimate of the graph, 
we prove the following. Note that $C^{1,\frac{1+\a}{2}}$ estimate has been established 
in \cite{Kasai-Tonegawa} and it will be integrated at the end. Some technical lemma 
concerning the change of second derivatives under orthogonal rotations is relegated to Section 8.
\begin{thm}
Corresponding to $n,\, k$, and $0<\alpha<1$ there exist $0<\varepsilon_3<1$ and $1<c_{13}<\infty$ 
with the following property. For $0<R<\infty$ suppose $\{V_t\}_{-R^2<t<R^2}$ and $\{u(\cdot,t)\}_{-R^2<
t<R^2}$, where $V_t=|M_t|$ with $M_t={\rm graph}\, f(\cdot,t)$, 
satisfy (B1)-(B4) on $C(T,R)\times (-R^2,R^2)$. Assume
\begin{equation}
\|\nabla f\|_{\frac{1+\alpha}{2}}:=\| \nabla f\|_{0}+R^{\frac{1+\alpha}{2}}[\nabla f]_{\frac{1+\alpha}{2}}\leq \varepsilon_3,
\label{s1}
\end{equation}
\begin{equation}
\|u\|_{\alpha}:=R\|u\|_{0}+R^{1+\alpha}[u]_{\alpha}\leq \varepsilon_3,
\label{s2}
\end{equation}
and assume that for some $g\in \F$ with 
\begin{equation}
\|\nabla g\|_0+R\|\nabla^2 g\|_0+R\|\frac{\p g}{\p t}\|_0\leq \varepsilon_3,
\label{s2.5}
\end{equation}
we have
\begin{equation}
\mu:=\Big(R^{-(k+4)}\int_{-R^2}^{R^2}\int_{C(T,R)}Q_{g}(x,t)\, d\|V_t\|(x)dt\Big)^{\frac12}\leq \varepsilon_3.
\label{s3}
\end{equation}
Then on $B_{R/2}^k\times(-R^2/4,R^2/4)$, 
$f$ is differentiable w.r.t. $(x,t)$ and $\nabla f$ is differentiable w.r.t. $x$,
and we have
\begin{equation}
R\Big\|\nabla^2 (f-g),\,\,\frac{\p (f-g)}{\p t}\Big\|_{0}+R^{1+\alpha}\Big[\nabla^2 f,\,\, \frac{\p f}{\p t}\Big]_{\alpha}
\leq c_{13}\max\{\mu,\|u\|_{\alpha},\|\nabla f\|_{\frac{1+\alpha}{2}}\}
\label{s4}
\end{equation}
where the (semi-)norms on the left-hand side of \eqref{s4} are over the domain $B_{R/2}^k$$\times$$(-R^2/4,R^2/4)$. 
Moreover, the normal velocity vector of $M_t$ is equal to $h(|M_t|, x)+u(x,t)^{\perp}$
at each point $x\in M_t \cap C(T,R/2)$ for $t\in (-R^2/4,R^2/4)$. 
\label{mthm}
\end{thm}
{\it Proof}. 
We may assume that $R=1$ after a change of variables. 
For any point $\tilde{x}\in M_{\tilde{t}} \cap C(T,1/2)$, $\tilde{t}\in (-1/2,1/2)$, there is a change of
variables by Lemma \ref{coor} so that the new graph function $\tilde{f}$ has $\nabla \tilde{f}(0,0)=0$ and 
$\tilde{u}(0,0)=0$ (where $(0,0)$ corresponds to $(\tilde{x},\tilde{t})$ before). Let $\tilde{T}$ be $\R^k\times\{0\}$
in this coordinate system which is also the tangent space to the graph of $\tilde{f}$ at $(0,0)$. 
 We will apply Proposition \ref{longprop} to $\tilde{f}$. To do so, we need the 
initial approximation function in $\F$. Consider the ${\rm graph}\, g$ after the
change of variables and let $\tilde{g}$ be the function defined on $\tilde{T}\times\R$ so that
${\rm graph}\, g={\rm graph}\, \tilde{g}$. Note that $\tilde{g}$ in general may not belong to $\F$.
Thus we do the following. In doing the above change of variables, choose a particular coordinate system so that it is
obtained first by the parallel translation $(\tilde{x},\tilde{t})\longmapsto (0,0)$, then by the change of variables 
$(x,t)\longmapsto (x-tu(\tilde{x},\tilde{t}),t)$, and 
an orthogonal rotation $A$ with $|I-A|\leq c(n,k)|\nabla f(\tilde{x},\tilde{t})|$
so that the image of the tangent space at $\tilde{x}$ of $M_{\tilde{t}}$ under $A$ is $\R^k\times\{0\}$. 
Define a polynomial function $\hat{g}$ with precisely the same first and second derivatives as $g$, that is,
if $g(x,t)=c+bt+\sum_{i=1}^k a_i x_i+1/2 \sum_{i,j=1}^k a_{ij}x_i x_j$ in the original coordinate
system, then we define $\hat{g}(x,t):=
bt+\sum_{i=1}^k a_i x_i+1/2\sum_{i,j=1}^k a_{ij}x_i x_j$. We emphasize to avoid 
any confusion that the variables $(x,t)$ for $\hat{g}$ is with respect to the new coordinate system. 
By definition, $\hat{g}\in \F$ and $\hat{g}(0,0)=\tilde{g}(0,0)$. Due to \eqref{nc1},
$|I-A|\leq c|\nabla f(\tilde{x},\tilde{t})|$ and
similar computations for the first derivatives, one has for some constant $c=c(n,k)$
\begin{equation}
\sup_{B_2^k \times(-1,1)}|\hat{g}-\tilde{g}|\leq c(|\nabla f(\tilde{x},\tilde{t})|+|u(\tilde{x},\tilde{t})|).
\label{s4.1}
\end{equation}
We then define
\begin{equation}
Q_{\hat{g}}(x,t):=\frac12\sum_{l=k+1}^n (x_l-\hat{g}_l(\tilde{T}(x),t))^2
\label{s4.2}
\end{equation}
where $\hat{g}=(\hat{g}_{l+1},\cdots,\hat{g}_n)$ and $x\in \R^n$ and similarly for
$Q_{\tilde{g}}$. By \eqref{s4.1}, we have for $x\in M_t$ with $t\in (-1,1)$ 
\begin{equation}
Q_{\hat{g}}\leq 2Q_{\tilde{g}}+c(n,k)(|\nabla f(\tilde{x},\tilde{t})|^2
+|u(\tilde{x},\tilde{t})|^2).
\label{s4.3}
\end{equation}
The difference between $Q_{\tilde{g}}$ and $Q_g$ on $M_t$ is that the former 
measure the $|\tilde{f}-\tilde{g}|^2/2$ while the latter measures 
$|f-g|^2/2$. The translation by $tu(\tilde{x},\tilde{t})$ does not affect the 
values of $Q_{\tilde{g}}$. Then a simple computation shows
\begin{equation}
Q_{\tilde{g}}\leq 2 Q_{g}.
\label{s4.4}
\end{equation}
Thus we have
\begin{equation}
\Big(\int_{-1/4}^{1/4}\int_{C(\tilde{T},1/4)}Q_{\hat{g}}(x,t)\, d\|V_t\|dt\Big)^{\frac12}
\leq 2 \mu+c(n,k)(\|\nabla f\|_0+\|u\|_0),
\label{ss4}
\end{equation}
where $V_t=|M_t|$ on the left-hand side is understood to be the one after the 
change of variables. Now we are in the position to apply Proposition \ref{longprop} 
for sufficiently small $\varepsilon_3$ which is determined by $\varepsilon_2$ and
\eqref{ss4}.
This proves that $\nabla \tilde{f}$ is differentiable w.r.t$.$ $x$ 
and $\tilde{f}$ is differentiable w.r.t$.$ $(x,t)$ at $(0,0)$. It is geometrically obvious that $f$ is then 
differentiable at $(\tilde{x},\tilde{t})$. It requires some calculations to prove that $\nabla f$ is 
differentiable w.r.t$.$ $x$ via computations as in Lemma \ref{changesec} but we omit the details. 
Moreover, since $g^{(0)}\in \F$, \eqref{cb2} proves that $\p \tilde{f}/\p t=\Delta \tilde{f}$ at $(0,0)$ for each component.
Since $\nabla \tilde{f}(0,0)=0$, this proves that the normal velocity is equal to the mean curvature at $(0,0)$.
Since the coordinate is `moving' with speed $u(\tilde{x},\tilde{t})$, we proved that the normal velocity is equal to the sum of the
mean curvature and $u^{\perp}$ in the original coordinate system. The supremum estimates for
$\nabla^2 (\tilde{f}-\hat{g})$ and $\p (\tilde{f}-\hat{g})/\p t$ follows from \eqref{cb3}.
This in turns 
gives
\begin{equation}
\sup_{B_{1/2}^k\times(-1/2,1/2)}\Big|\nabla^2(f-g),\,\frac{\p( f-g)}{\p t}\Big|\leq c_{14}\max\{\mu,\|u\|_{\a},\|\nabla f\|_{\frac{(1+\a)}{2}}\}
\label{s4.5}
\end{equation}
via \eqref{nc1} and estimates on the difference between $\nabla ^2\hat{g}$ and $\nabla^2\tilde{g}$,
which can be bounded by $c(n,k)|\nabla f|$.
Finally we need to prove the $\alpha$-H\"{o}lder norm estimate of \eqref{s4}. 
For $i=1,\, 2$, let $\tilde{x}_i\in M_{\tilde{t}_i} \cap C(T,1/2)$, $\tilde{t}_i\in (-1/2,1/2)$ be 
any two points with $(\tilde{x}_1,\tilde{t}_1)\neq (\tilde{x}_2,\tilde{t}_2)$. Without loss of
generality we assume
\begin{equation}
|\tilde{x}_1
-\tilde{x}_2|<1/10,\hspace{.3cm}\mbox{and}\hspace{.3cm}
0<\tilde{t}:=\tilde{t}_2-\tilde{t}_1<1/100.
\label{s5}
\end{equation}
After a change of variables 
as before, so that $(0,0)$
and $(\tilde{x},\tilde{t})$ in the new coordinate system correspond to 
$(\tilde{x}_1,\tilde{t}_1)$ and $(\tilde{x}_2,\tilde{t}_2)$, 
respectively, we may have $\nabla \tilde{f}(0,0)=0$ and
$\tilde{u}(0,0)=0$. Denote the tangent space to the graph $\tilde{f}$ at the origin by $\tilde{T}$. 
Restricting $\varepsilon_3$ further if necessary, 
by the first part of the proof and by Proposition \ref{longprop}, there exist $g^{(0)},\,
g^{(r)}\in {\mathcal F}$ for $0<r<1/4$ with \eqref{cb2}, \eqref{cb3} and \eqref{cb4} where $f$, $u$ and $T$
in those statements are replaced by $\tilde{f}$, $\tilde{u}$, $\tilde{T}$ with $R=1/4$. 
Corresponding to $(\tilde{x},\tilde{t})$, fix $\tilde{r}:=2\max\{|\tilde{x}|,|\tilde{t}|^{1/2}\}$
and consider $g^{(\tilde{r})}$. For later use, define
\begin{equation}
\tilde{a}:=(\tilde{a}_{ij})_{1\leq i,j\leq k}:=\big(\frac{\p^2 g^{(\tilde{r})}}{\p x_i\p x_j}\big)_{1\leq i,j\leq k},
\hspace{.5cm} \tilde{b}:=\frac{\p g^{(\tilde{r})}}{\p t}.
\label{s5.5}
\end{equation}
Note that (recall $\nabla g$ is independent of $t$ for $g\in \F$)
\begin{equation}
\begin{split}
|\nabla g^{(\tilde{r})}(\tilde{T}(\tilde{x}))| &\leq |\nabla g^{(\tilde{r})}(0)|+\tilde{r}
|\tilde{a}|\\ &\leq 2c_9 (4\tilde{r})^{1+\a}\max\{\mu,c_3[\tilde{u}]_{\a},c_3[\nabla \tilde{f}
]^2_{\frac{(1+\a)}{2}}\} +\tilde{r}\|\nabla^2 f\|_0 \\
&\leq c_{15}\tilde{r}\max\{\mu,\|u\|_{\a},\|\nabla f\|_{\frac{(1+\a)}{2}}\}
\end{split}
\label{s6}
\end{equation}
by \eqref{cb4}, the triangle inequalities and \eqref{s4.5}. Regarding the ${\rm graph}\, g^{(\tilde{r})}$ as a smooth
$k$-dimensional manifold in $\R^n$, let $\hat{T}\in {\bf G}(n,k)$ be the tangent space over
$(\tilde{x},\tilde{t})$ and let $\hat{g}$ be the graph representation over $\hat{T}$, that is, 
${\rm graph}\, \hat{g}={\rm graph}\, g^{(\tilde{r})}$. We introduce yet another new coordinate system
so that $\hat{T}=\R^k\times\{0\}$ and $(0,0)$ corresponds to $(\tilde{x},\tilde{t})$. We may take such new coordinate
system so that the new one is obtained by a parallel translation and an orthogonal rotation $A$ with
$|I-A|=O(|\nabla g^{(\tilde{r})}(\tilde{x})|)$. By \eqref{nc1} and \eqref{s6}, we have
\begin{equation}
\begin{split}
\sup_{(x,t)\in C(\hat{T},1/4)\times(-1/4,1/4)} & \big|\nabla^2 \hat{g}(x,t)-\tilde{a}\big|\leq
c(n,k)|\nabla g^{(\tilde{r})}(\tilde{x})||\tilde{a}| \\
&\leq c_{16}\tilde{r}\max\{\mu,\|u\|_{\a},\|\nabla f\|_{\frac{(1+\a)}{2}}\}.
\end{split}
\label{s7}
\end{equation}
Similar computations show
\begin{equation}
\sup_{(x,t)\in C(\hat{T},1/4)\times(-1/4,1/4)}\big|\frac{\p \hat{g}}{\p t}(x,t)-\tilde{b}\big|
\leq c_{16}\tilde{r}\max\{\mu,\|u\|_{\a},\|\nabla f\|_{\frac{(1+\a)}{2}}\}.
\label{s8}
\end{equation}
Now we define a function $\hat{g}^{(\tilde{r})}\in \F$ which is defined relative to $\hat{T}$ by
\begin{equation}
\hat{g}^{(\tilde{r})}(x,t):=\hat{g}(0,0)+\tilde{b}t+\frac12 \sum_{i,j=1}^k \tilde{a}_{ij}x_i x_j.
\label{s9}
\end{equation}
Since $g^{(\tilde{r})}\in \F$ and by \eqref{s5.5}, we have $\hat{g}^{(\tilde{r})}\in\F$. 
Moreover by the Taylor expansion and \eqref{s7}-\eqref{s9}, we have
\begin{equation}
\begin{split}
\sup_{B_{\tilde{r}}^k\times(-\tilde{r}^2,\tilde{r}^2)} & |\hat{g}-\hat{g}^{(\tilde{r})}|
\leq c(k)\tilde{r}^2\sup_{B_{\tilde{r}}^k \times(-\tilde{r}^2,\tilde{r}^2)}
\big(|\nabla^2 (\hat{g}- \hat{g}^{(\tilde{r})})|+\big|\frac{\p}{\p t}(\hat{g}-
\hat{g}^{(\tilde{r})})\big|\big) \\
& \leq 2c(k)c_{16}\tilde{r}^3 \max\{\mu,\|u\|_{\a},\|\nabla f\|_{\frac{1+\a}{2}}\}.
\end{split}
\label{s10}
\end{equation}
By 
$Q_{\hat{g}}\leq 2 Q_{g^{(\tilde{r})}}$ and \eqref{s10}, we have
\begin{equation}
\begin{split}\Big(\tilde{r}^{-k-4}\int_{-\tilde{r}^2/4}^{\tilde{r}^2/4}\int_{C(\hat{T},\tilde{r}/2)}
Q_{\hat{g}^{(\tilde{r})}}\, d\|V_t\|dt\Big)^{\frac12}
& \leq 2\Big(\tilde{r}^{-k-4}\int_{-\tilde{r}^2/4}^{\tilde{r}^2/4}\int_{C(\hat{T},\tilde{r}/2)}
Q_{g^{(\tilde{r})}}\, d\|V_t\|dt\Big)^{\frac12}\\
& +c(k)c_{16}^2\tilde{r}^2 \max\{\mu,\|u\|_{\a},\|\nabla f\|_{\frac{1+\a}{2}}
\} \end{split}
\label{s11}
\end{equation}
By \eqref{cb4}, the first term on the right-hand side is bounded
by $c_9 \tilde{r}^{1+\a}\max\{\mu,c_3 [u]_{\a},c_3 [\nabla f]^2_{\frac{1+\a}{2}}\}$.
At this point, we apply the first part of the present proof to conclude that the difference
between the second derivatives at $(0,0)$ and those of $\hat{g}^{(\tilde{r})}$ may be
bounded by a suitable constant multiples of 
$\tilde{r}^{\a}\max\{\mu, \|u\|_{\a}, \|\nabla f\|_{\frac{1+\a}{2}}\}$. The same holds for
time derivative. This proves the desired $\a$-H\"{o}lder estimate of \eqref{s4}.
\hfill{$\Box$}

Now we are in the position to prove our main Theorem \ref{mainreg2} and Theorem \ref{mainreg1}.
\newline
{\it Proof of Theorem \ref{mainreg2}}. 
As usual we may assume $R=1$. 
We apply \cite[Theorem 8.7]{Kasai-Tonegawa} first. To do so, we need to check the assumptions
(A1)-(A4) of \cite[Section 3.1]{Kasai-Tonegawa} are satisfied. Fix $p>k$ and $q>2$  as any large enough numbers
so that $\varsigma:=1-k/p-2/q>\frac{1+\a}{2}$. 
Since our $u$ is H\"{o}lder continuous, we have $\|u\|_{L^{p,q}}\leq \|u\|_0
\leq \|u\|_{\a}$ trivially. The upper bound (A2) can be proved via an argument in 
\cite[Proposition 6.2]{Kasai-Tonegawa}, or more specifically, one can show (with the
notations there)
\begin{equation}
\left.\int_{B_1}\hat{\rho}(\cdot,t)\, d\|V_t\|\right|_{t=t_1}^{t_2}\leq \int_{t_1}^{t_2}\left(\|u\|_0^2\int_{B_1}\hat{\rho}(\cdot,t)\, d\|V_t\|+c_1\omega_k E_0\right)\, dt
\label{ada1}
\end{equation}
using \eqref{acon0}. Solve a differential inequality for $\int_{B_1}\hat{\rho}\, d\|V_t\|$ using \eqref{ada1}. Moving 
around the location of pole, we obtain a uniform estimate (A2) in the interior. 
Thus, corresponding to the listed relevant constants, we have an interior
$C^{1,\varsigma}$ estimate for ${\rm spt}\, \|V_t\|$, i.e., we can represent ${\rm spt}\, \|V_t\|$
as a graph $f(\cdot,t)$ with the desired estimates for $\|f\|_0+\|\nabla f\|_{\frac{1+\a}{2}}$. 
Then use Theorem \ref{mthm} to obtain the second order derivatives estimates in a smaller
region, where we use $g=0$ for the initial approximation. 
Note that $\|\nabla f\|_{\frac{1+\a}{2}}$ on the right-hand side of 
\eqref{s4} is already estimated in terms
of $\mu$ and $\|u\|_{L^{p,q}}$. By choosing sufficiently small $\varepsilon_0>0$, this proves the desired conclusion.
\hfill{$\Box$}
\newline
{\it Proof of Theorem \ref{mainreg1}}.
Set $p,\, q$ as above. 
By the same reason as above, we have all the conditions (A1)-(A4) of \cite{Kasai-Tonegawa} satisfied.
Thus \cite[Theorem 3.2]{Kasai-Tonegawa} shows a.e$.$ $C^{1,\varsigma}$ regularity in space-time. 
Then Theorem \ref{mthm} shows $C^{2,\a}$ regularity there as well. 
\hfill{$\Box$}

\section{Brakke's MCF in submanifold}
It may be worthwhile to comment on some consequences of our main theorem
in the case that the ambient space $\R^n$ is replaced by a submanifold. 
Such situation naturally arises when we consider a MCF in general Riemannian manifold via
Nash's isometric imbedding theorem. For $\bar{k}\in \N$ with $1\leq k<\bar{k}\leq n$,
suppose we have a $C^{\infty}$ 
$\bar{k}$-dimensional submanifold $N$ in an open set $U\subset \R^n$ and a family of $k$-varifolds which 
is Brakke's MCF in $N$ in an appropriate weak sense.  
For the precise definition, we need to have a few preliminaries. 
We define the second fundamental form of $N$ at $x\in N$ 
to be the bilinear form ${\bf B}_x\,:\, {\rm Tan}_x N\times {\rm Tan}_x N\longmapsto ({\rm Tan}_x N)^{\perp}$ such that
\begin{equation}
{\bf B}_x(v_1,v_2):=-\sum_{i=1}^{n-\bar{k}}( v_1\cdot \nabla_{v_2}\tau_i)\tau_i\big|_x,\hspace{.5cm}
v_1,\,v_2\in {\rm Tan}_x N.
\label{r1}
\end{equation}
Here $\tau_1,\cdots,\tau_{n-\bar{k}}$ are locally defined vector fields which are orthonormal
and which satisfy $\tau_i (y)\in ({\rm Tan}_y N)^{\perp}$ on some neighborhood of $x$.
Next, for $x\in N$ and $S\in {\bf G}(n,k)$ with $S\subset {\rm Tan}_x N$, define
\begin{equation}
H_{N}(x,S):=\sum_{i=1}^k {\bf B}_x(v_i,v_i)\in ({\rm Tan}_x N)^{\perp},
\label{r2}
\end{equation}
where $v_1,\cdots,v_k$ is an orthonormal basis of $S$. $H_N(x,S)$ is well-defined independent 
of the choice of the orthonormal basis. 
Though it is simple, we record
the following
\begin{lemma}
Suppose $V\in {\bf IV}_k(U)$ satisfies ${\rm spt}\, \|V\|\subset N$ and has a generalized mean curvature
$h(V,\cdot)$ in $U$. Let $M\subset U$ be a countably $k$-rectifiable set such that $V=\theta|M|$
with some integer multiplicity function $\theta$. Then we have
\begin{equation}
h(V,x)-H_N(x,{\rm Tan}_x M)\in {\rm Tan}_x N
\label{r3}
\end{equation}
for $\H^k$ a.e$.$ on $M$. Here ${\rm Tan}_x M$ is the approximate tangent space of $M$ at $x$.
\label{pep}
\end{lemma}
{\it Proof}. It suffices to prove that 
\begin{equation}
\int_{U} (h(V,x)-H_N(x,S))\cdot f\, dV(x,S)=0
\label{r4}
\end{equation}
for all $f\in C^1_c(U;\R^n)$ with $f(x)\in ({\rm Tan}_x N)^{\perp}$. 
Let $\tau_1,\cdots,\tau_{n-\bar{k}}$ be a set of 
locally defined orthonormal vector fields which form a basis for $({\rm Tan}_y N)^{\perp}$
on $N$. Since the integration is over $M\subset N$, note that the values of $f$
outside of N do not matter. Thus without loss of generality we may express 
$f=\sum_{i=1}^{n-\bar{k}} f_i \tau_i$. Then by \eqref{first} we have
\begin{equation}
\int_U h(V,\cdot)\cdot f\, dV(\cdot,S)= -\sum_{i=1}^{n-\bar{k}} \int_U S\cdot 
\nabla (f_i\tau_i)\, dV(\cdot,S)=-\sum_{i=1}^{n-\bar{k}}\int_U
f_i S\cdot \nabla \tau_i\, dV(\cdot,S)
\label{r5}
\end{equation}
where we used $S\cdot \tau_i=0$ for $V$ a.e$.$ since $S={\rm Tan}_x M\subset
{\rm Tan}_x N$. On the other hand, by \eqref{r1} and \eqref{r2}, we see that $\int_U H_N(\cdot, S)
\cdot f\, dV(\cdot,S)$ is equal to the right-hand side of \eqref{r5}. This proves \eqref{r4}.
\hfill{$\Box$}
\begin{rem}
We should point out that $V$ being integral is not essential, and that it suffices for example to have
rectifiable $V$ with its approximate tangent space in ${\rm Tan}_x N$ a.e$.$ for Lemma \ref{pep}.
\end{rem}
Lemma \ref{pep} shows that for $V=\theta|M|$, we have a decomposition 
$h(V,x)=h(V,x)^{\top}+h(V,x)^{\perp}=h(V,x)^{\top}+H_N(x,{\rm Tan}_x M)\in {\rm Tan}_x N\oplus({\rm Tan}_x N)^{\perp}$.
Furthermore, due to the perpendicularity of the mean curvature vector \eqref{perpH}, 
we have $h(V,x)^{\top}\in ({\rm Tan}_x M)^{\perp}\cap {\rm Tan}_x N$ for
$\H^k$ a.e$.$ on $M$. The vector $h(V,\cdot)^{\top}$ may be considered as an intrinsic 
mean curvature vector with respect to $N$ and it is natural to define the mean curvature flow
whose velocity is equal to $h(V,\cdot)^{\top}$ as follows.
\begin{define}
 For $\Lambda\leq \infty$, a family of $k$-varifolds $\{V_t\}_{0\leq t<\Lambda}$ in $U\subset
 \R^n$ is (unit dentisy) Brakke's MCF in a smooth $\bar{k}$-dimensional submanifold
 $N\subset U$ if
 the followings are satisfied.
 \newline
 {\bf (C0)} For all $t\in [0,\Lambda)$, ${\rm spt}\, \|V_t\|\subset N$.
 \newline
 {\bf (C1)} For a.e$.$ $t\in [0,\Lambda)$, $V_t$ is a unit density $k$-varifold.
 \newline
 {\bf (C2)} For $\tilde{U}\subset\subset U$ and $(t_1,t_2)\subset\subset (0,\Lambda)$,
 \begin{equation}
 \sup_{t_1\leq t\leq t_2}\|V_t\|(\tilde{U})<\infty.
 \label{r6}
 \end{equation}
 {\bf (C3)} For all $\phi\in C^1 (N\times [0,\Lambda);\R^+)$ with $\phi(\cdot,t)\in C_c^1(N)$ and 
 $0\leq t_1<t_2<\Lambda$, we have
 \begin{equation}
 \begin{split}
 \|V_{t_2}\|&(\phi(\cdot,t_2))-\|V_{t_1}\|(\phi(\cdot,t_1)) \\ &\leq \int_{t_1}^{t_2}\int_{G_k(U)} (\nabla_N \phi-\phi h(V_t,\cdot))
 \cdot(h(V_t,\cdot)-H_N(\cdot,S))+\frac{\p\phi}{\p t}(\cdot, t)\, dV_t(\cdot,S)dt,
 \end{split}
 \label{r7}
\end{equation}
where $\nabla_N\phi$ is the tangential derivative of $\phi$ on $N$.
\end{define}
\begin{rem} We also assume that $h(V_t,\cdot)$ exists for a.e$.$ $t$ and locally $L^2$ integrable with respect to
$d\|V_t\|\,dt$. 
By (C1) and Lemma \ref{pep}, for a.e$.$ $t$, we may replace
both the first $h(V,\cdot)$ and $(h(V,\cdot)-H_N(\cdot,S))$ of \eqref{r7} by $h(V,\cdot)^{\top}$ without changing the definition.
We may also ask (C3) to hold for $\phi$ defined on $U$ and for $\nabla\phi$ instead of $\nabla_N\phi$ due to Lemma \ref{pep}. In sum, we may equivalently assume the following. 
\newline
{\bf (C3)'} For all $\phi\in C^1(U\times[0,\Lambda);\R^+)$ with $\phi(\cdot,t)\in C^1_c(U)$ and $0\leq t_1<t_2<\Lambda$,
we have
\begin{equation}
\begin{split}
\|V_{t_2}\|&(\phi(\cdot,t_2))-\|V_{t_1}\|(\phi(\cdot,t_1)) \\ &\leq \int_{t_1}^{t_2}\int_{G_k(U)}
(\nabla\phi-\phi h(V_t,\cdot)^{\top})\cdot h(V_t,\cdot)^{\top}+\frac{\p \phi}{\p t}(\cdot,t)\, dV_t(\cdot,S)dt.
\end{split}
\label{r8}
\end{equation} 
\end{rem}
Now let us discuss what can be said under the assumptions (C0)-(C3).
Since $H_N(\cdot,S)$ is locally a bounded function with $H_N(x,S)\in ({\rm Tan}_x N)^{\perp}$, we may
regard $H_N$ as $u^{\perp}$ in \cite[Theorem 3.2]{Kasai-Tonegawa} for any large $p$ and $q$. Thus we may
conclude that $M_t:={\rm spt}\, \|V_t\|$ is a $C^{1,\varsigma}$ graph for a.e$.$ in space-time. This in turn shows
that $H_N(x,{\rm Tan}_x M_t)$ is $\varsigma$-H\"{o}lder continuous since it involves the first derivatives of the
graph. This will lead us to the setting of the present paper, which shows partial $C^{2,\alpha}$ regularity with
motion law `velocity $= h(V_t,\cdot)^{\top}$' being satisfied classically. Then the standard parabolic regularity
theory shows partial $C^{\infty}$ regularity. Thus we proved
that {\it any unit density Brakke's MCF in submanifold is necessarily a.e$.$ smooth}, the 
meaning of a.e$.$ is stated rigorously in Section \ref{pr}. The corresponding statement for
general smooth Riemannian manifold setting also follows via Nash's imbedding theorem. 
\section{Appendix}
In this appendix we consider how the second derivatives change under the orthogonal change of variables.
\begin{lemma} There exist $\b=\b(n,k)$ and $c=c(n,k)$ with the following property. 
Suppose that $A=(a_{ij})_{1\leq i,j\leq n}$ is an orthogonal matrix with
\begin{equation}
|I-A|\leq \b.
\label{ncsmall}
\end{equation}
Suppose that two coordinate systems $x=(x_1,\cdots, x_n)$ and
$\tilde{x}=(\tilde{x}_1,\cdots,\tilde{x}_n)$ are related by $\tilde{x}^t=Ax^t$. 
Suppose that a $k$-dimensional manifold $M$  
in $\R^n$ is represented in the $x$ and $\tilde{x}$
coordinate systems as $x_{j}=f_{j}(x_1,\cdots,x_k)$ for $j=k+1,\cdots,n$
and $\tilde{x}_{j}=\tilde{f}_{j}(\tilde{x}_1,
\cdots,\tilde{x}_k)$ for $j=k+1,\cdots,n$, respectively. Assume that $f_{j}$ and
$\tilde{f}_j$ are differentiable for $j=k+1,\cdots,n$. Further assume that
\begin{equation}
|\nabla f|,\, |\nabla \tilde{f}|\leq 1
\label{nc0}
\end{equation}
on the domain of definitions of $f$ and $\tilde{f}$. 
Let two points in $M$ be expressed in $x$ coordinate system as 
$(x_1^{(i)},\cdots,x_k^{(i)},f_{k+1}(x_1^{(i)},\cdots,x_k^{(i)}),\cdots,
f_{n}(x_1^{(i)},\cdots,x_k^{(i)}))$ for $i=1,2$ and in $\tilde{x}$ coordinate
system as $(\tilde{x}_1^{(i)},\cdots,\tilde{x}_k^{(i)},
\tilde{f}_{k+1}(\tilde{x}_1^{(i)},\cdots,\tilde{x}_k^{(i)}),\cdots,
\tilde{f}_{n}(\tilde{x}_1^{(i)},\cdots,\tilde{x}_k^{(i)}))$ for $i=1,2$, respectively. Then
writing $x^{(i)}=(x_1^{(i)},\cdots,x_k^{(i)})$ and $\tilde{x}^{(i)}=(\tilde{x}_1^{(i)},
\cdots,\tilde{x}_k^{(i)})$, we have
\begin{equation}
|\nabla f(x^{(1)})-\nabla f(x^{(2)})|\leq c(n,k) |\nabla \tilde{f}(\tilde{x}^{(1)})-\nabla \tilde{f}
(\tilde{x}^{(2)})|.
\label{nc3}
\end{equation}
Furthermore, assume that $f_j$ and $\tilde{f}_j$ are twice differentiable for $j=k+1,\cdots,
n$. Then we have
\begin{equation}
|\nabla^2 f(x^{(i)})-\nabla^2 \tilde{f}(\tilde{x}^{(i)})|\leq c(n,k)|I-A||\nabla^2 \tilde{f}(\tilde{x}^{(i)})|,
\label{nc1}
\end{equation}
\begin{equation}
\begin{split}
|\nabla^2 f(x^{(1)})-\nabla^2 f(x^{(2)})| &
\leq c(n,k)|I-A||\nabla \tilde{f}(\tilde{x}^{(1)})-\nabla \tilde{f}(\tilde{x}^{(2)})|\max_{i=1,2}
|\nabla^2 \tilde{f}(\tilde{x}^{(i)})| \\
&+c(n,k)|\nabla^2 \tilde{f}(\tilde{x}^{(1)})-\nabla^2 
\tilde{f}(\tilde{x}^{(2)})|.
\end{split}
\label{nc2}
\end{equation}
\label{changesec}
\end{lemma}
{\it Proof}.
For the moment we drop the upper subscript $(i)$ for simplicity. Since two coordinate systems are
related by $\tilde{x}^t=Ax^t$, we have
\begin{equation}
(\tilde{x}_1,\cdots,\tilde{x}_k, \tilde{f}_{k+1},\cdots,
\tilde{f}_n)^t=A (x_1,\cdots,x_k, f_{k+1},
\cdots,f_n)^t.
\label{n1}
\end{equation}
By \eqref{n1}, one obtains the following identity for each $m=k+1,\cdots,n$, 
\begin{equation}
\sum_{l\leq k} a_{ml}x_l+\sum_{l>k}a_{ml} f_l=
\tilde{f}_m\big(\sum_{l\leq k}a_{1l}x_l+\sum_{l>k}a_{1l} f_{l},
\cdots, \sum_{l\leq k}a_{kl}x_l+\sum_{l>k}a_{kl}f_{l}\big).
\label{n2}
\end{equation}
Differentiating \eqref{n2} with respect to $x_i$, $1\leq i\leq k$, we have (writing 
$f_{l,i}:=\frac{\p f_l}{\p x_i}$ and similarly for $\tilde{f}$)
\begin{equation}
a_{mi}+\sum_{l>k}a_{ml} f_{l,i}=\sum_{p\leq k}(a_{pi}+\sum_{l>k}
a_{pl}f_{l,i})\tilde{f}_{m,p}.
\label{n2.5}
\end{equation}
Assuming that $f_j$ and $\tilde{f}_j$ are twice differentiable, and differentiating \eqref{n2.5} with respect to $x_j$, $1\leq j\leq k$, we have
(writing $f_{l,ij}:=\frac{\p^2 f_l}{\p x_i x_j}$ and similarly for $\tilde{f}$)
\begin{equation}
\sum_{l>k}a_{ml}f_{l,ij}=\sum_{p,q\leq k}(a_{pi}+\sum_{l>k}
a_{pl}f_{l,i})(a_{qj}+\sum_{l'>k}a_{ql'}f_{l',j})\tilde{f}_{m,pq}+
\sum_{p\leq k}\tilde{f}_{m,p}\sum_{l>k}a_{pl}f_{l,ij}.
\label{n3}
\end{equation}
Moving the last term of \eqref{n2.5} and \eqref{n3} to the right-hand side, respectively, we obtain
\begin{equation}
\sum_{l>k}(a_{ml}-\sum_{p\leq k}\tilde{f}_{m,p} a_{pl})f_{l,i}=-a_{mi}
+\sum_{p\leq k}a_{pi}\tilde{f}_{m,p},
\label{n3.5}
\end{equation}
\begin{equation}
\sum_{l>k}(a_{ml}-\sum_{p\leq k}\tilde{f}_{m,p} a_{pl})f_{l,ij}=\sum_{p,q\leq k}
(a_{pi}+\sum_{l>k}a_{pl}f_{l,i})(a_{qj}+\sum_{l'>k}a_{ql'}f_{l',j})\tilde{f}_{m,pq}.
\label{n4}
\end{equation}
Define $(n-k)\times(n-k)$ matrix-valued function
$\tilde{A}=(\tilde{A}_{ij})_{1\leq i,j\leq n-k}$ whose $(i,j)$ component is 
$a_{i+k,\,j+k}-\sum_{p\leq k}\tilde{f}_{i+k,p} a_{p,\,j+k}$. By \eqref{ncsmall} and
\eqref{nc0}, if we restrict $\e$ sufficiently small, $\tilde{A}$ is 
invertible. Multiplying $\tilde{A}^{-1}$ from left to \eqref{n3.5} and \eqref{n4}, respectively, 
we obtain 
\begin{equation}
f_{m,i}=-\sum_{m'>k}(\tilde{A}^{-1})_{mm'}a_{m' i}
+\sum_{p\leq k,\,m'>k}(\tilde{A}^{-1})_{mm'}a_{pi}\tilde{f}_{m', p},
\label{n4.5}
\end{equation}
\begin{equation}
f_{m,ij}=\sum_{p,q\leq k}(a_{pi}+\sum_{l>k}a_{pl}f_{l,i})
(a_{qj}+\sum_{l'>k}a_{ql}f_{l,j}) \sum_{m'>k}(\tilde{A}^{-1})_{mm'}\tilde{f}_{m',pq}.
\label{n5}
\end{equation}
For \eqref{nc3}, it is not difficult to check (by the definition of inverse matrix)
that
\begin{equation}
|(\tilde{A}^{-1})_{mm'}(\tilde{x}^{(1)})-(\tilde{A}^{-1})_{mm'}(\tilde{x}^{(2)})|
\leq c(n,k)|I-A||\nabla \tilde{f}(\tilde{x}^{(1)})-\nabla \tilde{f}(\tilde{x}^{(2)})|.
\label{n5.5}
\end{equation}
Then \eqref{nc3} follows from \eqref{n4.5}, \eqref{n5.5} and the triangle inequality.
Note that we only need differentiability to obtain \eqref{nc3}. 
We next consider the difference between $f_{m,ij}$ and $\tilde{f}_{m,ij}$.
The sum of the right-hand side of \eqref{n5} is separated to $\sum_{(p,q)=(i,j)}\cdot+\sum_{(p,q)\neq (i,j)}\cdot
=: E_1+E_2$. Then
\begin{equation}
\begin{split}
|E_1-\tilde{f}_{m,ij}|&=
\big|(a_{ii}+\sum_{l>k}a_{il}f_{l,i})(a_{jj}+\sum_{l'>k}a_{jl'}f_{l',j})\sum_{m'>k}(\tilde{A}^{-1})_{mm'}\tilde{f}_{m',ij}
-\tilde{f}_{m,ij}\big|\\
& \leq |\tilde{f}_{m,ij}|\big|1-(a_{ii}+\sum_{l>k}a_{il}f_{l,i})(a_{jj}+\sum_{l'>k}a_{jl'}f_{l',j})(\tilde{A}^{-1})_{mm}|\\
&+|(a_{ii}+\sum_{l>k}a_{il}f_{l,i})(a_{jj}+\sum_{l'>k}a_{jl'}f_{l',j})\sum_{m'\neq m}|(\tilde{A}^{-1})_{mm'}\tilde{f}_{m',ij}|\\
&\leq c(n,k)|I-A|\sum_{m'>k}|\tilde{f}_{m',ij}|
\end{split}
\label{n6}
\end{equation}
since the off-diagonal elements of $A$ and $\tilde{A}^{-1}$ are bounded by $c(n,k)|A-I|$.
For $E_2$, since $p\neq i$ or $q\neq j$, one can check from \eqref{n5} that 
\begin{equation}
|E_2|\leq c(n,k)|I-A|\sum_{m'>k,\, p,q\leq k}|\tilde{f}_{m',pq}|.
\label{n7}
\end{equation}
\eqref{n5}-\eqref{n7} prove \eqref{nc1}. For \eqref{nc2}, we have for $p\leq k$
\begin{equation}
|\sum_{l>k} a_{pl}(f_{l,i}(x^{(1)})-f_{l,i}(x^{(2)}))|\leq c(n,k)|I-A||\nabla f(x^{(1)})
-\nabla f(x^{(2)})|.
\label{n8}
\end{equation} 
Then \eqref{n5}, \eqref{n5.5}, \eqref{n8} (with \eqref{nc3}) with
suitable triangle inequalities prove \eqref{nc2}.
\hfill{$\Box$}

\end{document}